\documentclass[bezier,12pt]{article}
\usepackage{amssymb}
\usepackage{latexsym}
\usepackage{color}
\usepackage{amsmath}
\usepackage{amsfonts}

\newtheorem{Theorem}{Theorem}[part]

\newtheorem{Proposition}{Proposition}[part]

\newtheorem{Lemma}{Lemma}[part]
\newtheorem{Corollary}{Corollary}[part]
\newtheorem{Remark}{Remark}[part]

\makeatletter \@addtoreset{equation}{section}

\@addtoreset{Definition}{section}

\@addtoreset{Theorem}{section}

\@addtoreset{Proposition}{section}

\@addtoreset{Property}{section}

\@addtoreset{Assumption}{section}

\@addtoreset{Corollary}{section}

\@addtoreset{Lemma}{section}

\@addtoreset{Remark}{section}

\@addtoreset{Example}{section}

\def \proof{{\noindent \bf Proof. }}
\def \ep{\hbox{ }\hfill$\Box$}
\def\reff#1{{\rm(\ref{#1})}}

\def\Ac{{\cal A}}

\def\eps{\varepsilon}

\def\Bc{{\cal B}}
\def\Dc{{\cal D}}

\def\Fc{{\cal F}}
\def\Gc{{\cal G}}
\def\Hc{{\cal H}}
\def\Jc{{\cal J}}

\def\Jc{{\cal J}}

\def\Pc{{\cal P}}

\def\Uc{{\cal U}}

\def\eps{\varepsilon}

\def\diag{\mbox{\rm diag}}

\def\no{\noindent}
\def\Pas{\mathbb{P}-\mbox{a.s.}}

\def\x{\times}

\def\05{\frac{1}{2}}
\def\-1{^{-1}}
\def\1{{1\hspace{-1mm}{\rm I}}}
\def\={\;=\;}
\def\.{\;.}
\def\vp{\varphi}
\title{  }
\author{ }

\def\be{\begin{eqnarray}}
\def\ee{\end{eqnarray}}

\def\lra{\longrightarrow}

\def\b*{\begin{eqnarray*}}
\def\e*{\end{eqnarray*}}

\def\And{\;\mbox{ and }\;}

\def \ep{\hbox{ }\hfill{ ${\cal t}$~\hspace{-5.5mm}~${\cal u}$   } }

\def\Esp#1{\mathbb{E}\left[#1\right]}

\def\Pro#1{\mathbb{P}\left[{#1}\right]}

\def\And{\mbox{ and } }
\def\pourtout{\mbox{ for all } }

\def \E{\mathbb{E}}
\def \F{\mathbb{F}}
\def \M{\mathbb{M}}
\def \N{\mathbb{N}}
\def \R{\mathbb{R}}

\def\P{\mathbb{P}}

\def\lb{\underline}

\def \F{I\!\!F}
\def \M{I\!\!M}
\def \N{I\!\!N}
\def \R{I\!\!R}

\def\Ac{{\cal A}}
\def\Bc{{\cal B}}

\def\Dc{{\cal D}}

\def\Fc{{\cal F}}
\def\Gc{{\cal G}}
\def\Hc{{\cal H}}
\def\Jc{{\cal J}}

\def\Pc{{\cal P}}

\def\Uc{{\cal U}}

\def\lra{\longrightarrow}

\def\-1{^{-1}}

\def\0.5{\frac{1}{2}}

\def\x{\times}
\def\no{\noindent}
\def\={\;=\;}
\def\.{\;.}

\def \proof{{\noindent \bf Proof. }}
\def\reff#1{{\rm(\ref{#1})}}
\def\eps{\varepsilon}
\def\vp{\varphi}

\def\diag{\mbox{\rm diag}}

\def\1{{\bf 1}}
\def\Tr{\mbox{\rm Tr}}
\def \ep{\hbox{ }\hfill{ ${\cal t}$~\hspace{-5.5mm}~${\cal u}$   } }

\title{Explicit characterization of the super-replication strategy in
financial markets with partial transaction costs}

\author{ Imen Bentahar  \\\small Universit\'e Dauphine, CEREMADE,
                        \\\small Paris, France
                        \\\small imen.ben\_tahar@ceremade.fr
\and Bruno Bouchard     \\\small Universit\'e Paris VI, LPMA,
                        \\\small and CREST
                        \\\small Paris, France
                        \\\small bouchard@ccr.jussieu.fr}

\begin{document}
\vspace{50pt} \maketitle \vspace{.5in}

\begin{abstract}
We consider a multivariate financial market with transaction
costs and study the problem of finding the minimal initial capital
needed to hedge, without risk, European-type contingent claims.
The model is similar to the one considered in  Bouchard and Touzi (2000), except
that some of the assets can be exchanged freely, i.e. without
paying transaction costs. In this context, we generalize the
result of the above paper and prove that the value of this stochastic
control problem is given by the cost of the cheapest hedging
strategy in which the number of non-freely exchangeable assets is
kept constant over time.
\end{abstract}

\vspace{10mm}

\noindent{\sl Kew words:} Non efficient transaction costs, hedging
options,   viscosity solutions.


\newpage

\section{Introduction}

Since the 90's, there has been many papers devoted to the proof
of the conjecture of Davis and Clark (1994)~: in the context of
the Black and Scholes model with proportional transaction costs,
the cheapest super-hedging strategy for a European call option is
just the price (up to initial transaction costs) of the
underlying asset. The first proofs of this result were obtained,
independently, by Soner, Shreve and Cvitani\'c (1995) and
Levental and Skorohod (1995). In a one-dimensional Markov
diffusion model, the result was extended
 by Cvitani\'c, Pham and Touzi (1999) for general contingent
claims. Their approach relies on the dual formulation of the super-replication cost
(see Jouini and Kallal 1995 and Cvitani\'c and Karatzas 1996).

The multivariate case was then considered by Bouchard and Touzi
(2000). In contrast to Cvitani\'c, Pham and Touzi (1999), they
did not use the dual formulation but  introduced a family of
fictitious markets without transaction costs but with modified
price processes evolving in the bid-ask spreads of the original
market. Then, they defined the associated super-hedging problems
and showed that they provide lower bounds for the original one. By
means of a direct dynamic programming principle for stochastic
targets problems, see e.g. Soner and Touzi (2002), they provided
a PDE characterization for the upper bound of these auxiliary
super-hedging prices. Using similar arguments as in Cvitani\'c,
Pham and Touzi (1999), they were then able to show that the
associated value function is concave in space and non-increasing
in time. This was enough to show that it corresponds to the price
of the cheapest buy-and-hold strategy in the original market. A
different proof relying on the dual formulation for multivariate
markets, see Kabanov (1999), was then proposed  by Bouchard
(2000).

It should be noticed that a crucial point of all this analysis is
that transaction costs are efficient, i.e. there is no couple of
freely exchangeable assets. In this paper, we propose a first
attempt to characterize the super-replication strategy in
financial markets with ``partial" transaction costs, where some
assets can be exchanged freely. As a first step, we follow the
approach of Bouchard and Touzi (2000). We introduce a family
of fictitious markets and provide a PDE characterization similar
to the one obtained in this paper. However, in our context, one
can only show that the corresponding value function is concave in
some directions (the ones where transaction costs are effective),
and this is not sufficient to provide a precise characterization
of the super-hedging strategy. With the help of a comparison
theorem for PDE's, we next obtain a new lower bound associated to
an auxiliary control problem, written in standard from. This
allows us to characterize the optimal hedging strategy~: it
consists in keeping constant the number of non-freely exchangeable
assets held in the portfolio and hedging the remaining part of the
claim by trading dynamically on the freely exchangeable ones.

The paper is organized as follows. After setting some notations
in Section \ref{sec Nota}, we describe the model and the
super-replication problem in Section \ref{sec model}. The main
result of the paper is stated in Section \ref{sec thmmain}. In
Section \ref{sec marche fictif}, we   introduce an auxiliary
super-hedging problem similar to the one considered in Bouchard
and Touzi (2000) and derive the PDE associated to the value
function. Further properties of this value function are obtained
in Section \ref{sec pbm auxiliaire}. The proof is concluded in
Section \ref{sec proof finale thmmain}.

\section{Notations}\label{sec Nota}

For the reader's convenience, we first introduce the main notations of this paper.

\no Given a vector $x$ $\in$ $\R^n$, its $i$-th component is
denoted by $x^i$. All elements of $\R^n$ are identified with
column vectors and the scalar product is denoted by $\cdot$.
$\M^{n,p}$ denotes the set of all real-valued matrices with $n$
rows and $p$ columns. Given a matrix $M$ $\in$ $\M^{n,p}$, we
denote by $M^{ij}$ the component corresponding to the $i$-th row
and the $j$-th column. $\M^{n,p}_+$ denotes the subset of
$\M^{n,p}$ whose elements have non-negative entries. If $n=p$, we
simply denote $\M^n$ and $\M^n_+$ for $\M^{n,n}$ and
$\M^{n,n}_+$. Since $\M^{n,p}$ can be identified with $\R^{np}$,
we define the norm on $\M^{n,p}$ as the norm of the associated
element of $\R^{np}$. In both cases, the norm is simply denoted
by $|\cdot|$. Transposition is denoted by $'$.  Given a square
matrix $M$ $\in$ $\M^{n}$, we denote by Tr$[M]$ $:=$
$\sum_{i=1}^nM^{ii}$ the associated trace. For $x$ $\in \R^p$ and
$\eta>0$, we denote by $B(x,\eta)$ the open ball with radius
$\eta$ centered in $x$, $\partial B(x,\eta)$ its boundary and
$\bar B(x,\eta)$ its closure.

\no Given n scalars $x_1,\ldots,x_n$,
we denote by Vect$[x_i]_{i\le n}$ the vector
of $\R^n$ defined by the components $x_1,\ldots,x_n$.
For all $x$ $\in$ $\R^n$, diag$[x]$ denotes the diagonal
matrix of $\M^n$ whose $i$-th diagonal element is $x^i$.

\no We denote by ${\bf 1}_i$ the vector of $\R^n$ defined by
${\bf 1}_i^j$ $=$ $1$ if $j$ $=$ $i$ and $0$ otherwise.

\no Given a smooth function $\varphi$ mapping $\R^n$ into $\R^p$,
we denote by $D_z \varphi$ the (partial) Jacobian matrix of
$\varphi$ with respect to its $z$ variable. In the case $p$ $=$
$1$, we denote by $D^2_{zs}\varphi$ the matrix defined as
$(D^2_{zs}\varphi)^{ij}$ $=$ $\partial^2 \varphi/\partial
z^i\partial s^j$. If $\vp$ depends only on $z$, we simply write
$D\vp$ and $D^2\vp$ in place of  $D_z\vp$ and $D^2_{zz}\vp$.

\no In this paper, we shall consider $\R^d-$valued variables,
$d\ge 1$, and we shall often write $d$ as the sum of two positive
integers $d_f+d_c$.  For $x$ $\in$ $\R^d$, we will then write $x$
as $x$ $=$ $(x^f,x^c)$ where $x^f$ (resp. $x^c$) stands for the
vector of $\R^{d_f}$ (resp. $\R^{d_c}$) formed by the $d_f$ first
(resp. $d_c$ last) components of $x$.

\no All inequalities involving random variables have to be
understood in the $\Pas$ sense.

\section{The model}\label{sec model}

Let $T$ be a finite time horizon and $(\Omega,\Fc,\P)$ be a
complete probability space supporting a $d$-dimensional Brownian
motion $\{W(t)$, $0\leq t\leq T\}$. We shall denote by $\F$ $=$
$\{\Fc_t$, $0\leq t\leq T\}$ the $\P$-augmentation of the
filtration generated by $W$.

\subsection{The financial market}

\no We consider a financial market which consists of one bank
account, with constant price process, normalized to unity, and $d$
risky assets $S:=\{S^1,\ldots,S^d\}'$. The price process $S$ $=$
$\{S(t)$, $0\leq t\leq T\}$ is an $\R^{d}$-valued stochastic
process defined by the following stochastic differential system
\begin{eqnarray}\label{def_s}
dS(t) &=& \mbox{\rm diag}[S(t)]\sigma (t,S(t))dW(t)  ,
\hspace{10mm} 0 \;<\; t \;\leq \;T \;.
\end{eqnarray}
Here $\sigma (.,.)$ is an $\M^{d}$-valued function. We shall
assume all over the paper that the function diag$[s]\sigma(t,s)$
satisfies the usual Lipschitz and linear growth conditions in
order for the process $S$ to be well-defined  and that
$\sigma(t,s)$ is invertible with $\sigma(t,s)^{-1}$ locally
bounded, for all $(t,s)\in [0,T]\times \R^d_+$.



\begin{Remark}
{\rm As usual, the assumption that the interest rate of the bank
account is zero could be easily dispensed with by discounting.
Also, there is no loss of generality in defining $S$  as a
martingale since we can always reduce the model to this context
by an appropriate change of measure (under mild conditions on the
initial coefficients). }
\end{Remark}

\no We write $d$ as $d=d_f+d_c$ with $d_f,\;d_c\ge 1$. The
subscript $f$ stands for ``free" (of costs) while $c$ means
``costs".  We assume that transactions in the market formed by
the first $d_f$ assets and the num\'eraire are free of costs. This
means that a portfolio process associated to a trading strategy
$\phi$ only based on the first $d_f$ assets can be written  in
the usual form $x+\int_0^\cdot \phi(t)\cdot dS^f(t)$, where, for
$z \in \R^d$, we recall that $z^f$ $=$ $(z^1,\ldots,z^{d_f})'$
and $z^c$ $=$ $(z^{d_f+1},\ldots,z^{d})'$. On the other hand, we
assume that any transaction involving the last $d_c$ assets is
subject to proportional costs.

\no Then, a trading strategy is described by a pair $(\phi,L)$
where $\phi$ is a $\R^{d_f}$-valued predictable process satisfying
$\int_0^T|\phi(t)|^2 dt <\infty$  and $L$ is an
$\M^{1+d_c}_+$-valued process   with initial value $L(0-)$ $=$
$0$, such that $L^{ij}$ is $\F-$adapted, right-continuous, and
nondecreasing for all $i,j$ $=$ $1,\ldots,1+d_c$.

\no To such a pair $(\phi,L)$   and $x$ $\in$ $\R^{1+d_c}$, we
then associate the wealth process $X_x^{\phi,L}$ defined as the
solution on $[0,T]$ of
    \b*
    X_{x}^{i}(0-) &=& x^{i} \;\quad\;\mbox{ for  } 1\le i\le 1+d_c\\
    dX_{x}^{1+i}(t) &=& X_{x}^{1+i}(t)\ \frac{dS^{d_f+i}(t)}{S^{d_f+i}(t)}
    +\sum_{j=1}^{1+d_c}\left[ dL^{j(1+i)}(t)-(1+\lambda ^{(1+i)j})dL^{(1+i)j}(t)\right]
    \\
    && \hspace{9cm}\;\;\mbox{ for } 1\le i\le d_c\;,\\
    dX_{x}^{1}(t) &=& \phi(t)\cdot   dS^{f}(t) +
    \sum_{j=1}^{1+d_c}\left[ dL^{j1}(t)-(1+\lambda^{1j})dL^{1j}(t)\right]
    \;.
    \e*
The first component $(X_x^{\phi,L})^1$ stands for the amount of
money which is invested on the market formed by the $d_f$ first
assets  and the num\'eraire. For $i>1$, $(X_x^{\phi,L})^i$ is the
amount invested in the $(d_f+i-1)$-th asset. The coefficient
$\lambda^{ij}$ are assumed to be non-negative and stands for the
transaction costs which are paid when transacting on one of the
$d_c$ last assets.

\no Observe that we can always assume that
    \b*
    (1+\lambda^{ij})&\le&
    (1+\lambda^{ik})(1+\lambda^{kj})\;\;,\;i,\;j,\;k\le 1+d_c\;,
    \e*
since otherwise it would be cheaper to transfer money from the
account $i$ to $j$ by passing through $k$ rather than directly.
Then, for any ``optimal" strategy the effective cost between $i$
and $j$ would be $(1+\lambda^{ik})(1+\lambda^{kj})-1$.


\subsection{The super-replication problem}

Following Kabanov (1999), we define the {\it solvency region}~:
\begin{eqnarray*}
K:= \left\{\; x\in \R^{1+d_c}~:~\exists\; a\in\M_{+}^{1+d_c},\;
                               x^{i}+\sum_{j=1}^{1+d_c}(a^{ji}-(1+\lambda^{ij})a^{ij})
                               \;\geq\; 0 \; \;\forall\; i\le 1+d_c \;\right\}\;.
\end{eqnarray*}
The elements of $K$ can be interpreted as the vectors of portfolio
holdings such that the no-bankruptcy condition is satisfied, i.e.
the liquidation value of the portfolio holdings $x,$ through some
convenient transfers, is nonnegative.

\smallskip

\no Clearly, the set $K$ is a closed convex cone containing the
origin. We then introduce the partial ordering $\succeq$ induced
by $K$~:
\[
\mbox{for all }x_{1},x_{2}\in \R^{1+d_c},\mbox{
\quad }{\normalsize x_{1}\succeq x_{2}}\mbox{ }{\normalsize \ }\mbox{%
{\normalsize if and only if }}{\normalsize \ x_{1}-x_{2}\in K \;.\
}
\]
A trading strategy $(\phi,L)$ is said to be {\it admissible} if
there is some some $c,\;\delta\in \R$ and $\delta_f$ in $\R^{d_f}$ such that
\begin{eqnarray}
\label{posi_x} X_{0}^{\phi,L}(t) &\succeq& - (c+\delta_f\cdot
S^f(t),\delta S^c(t))\; ,\;\; 0\leq t\leq T \;.
\end{eqnarray}
Observe that, if $X_{0}^{\phi,L}$ satisfies the above condition,
then, after possibly changing $(c,\delta,\delta_f)$, it holds for
$X_{x}^{\phi,L}$ too, $x\in \R^{1+d_c}$. We  denote by
${\cal A}$ the set of such trading strategies.
\\

\no A {\it contingent claim} is a $(1+d_c)$-dimensional
$\Fc_T$-measurable random variable $g(S(T))$. Here, $g$ maps
$\R^d_+$ into $\R^{1+d_c}$ and satisfies
    \begin{equation}
    g(s)\succeq -  (c+\delta_f\cdot s^f,\delta s^c)\quad \quad\mbox{for all }s\mbox{ in }\R^d_+ \; \label{TA}
    \end{equation}
for some   $c,\;\delta\in \R$ and $\delta_f$ in $\R^{d_f}$ .

\no In the rest of the paper, we shall identify a contingent claim
with its pay-off function $g.$ For  $i=2,\ldots,1+d_c$, the
random variable $g^i(S(T))$ represents a target position in the
asset $d_f-1+i$, while $g^1(S(T))$ represents a target position
in terms of the num\'eraire.
\\

\no The {\it super-replication problem} of the contingent claim
$g$ is then defined by
\begin{eqnarray*}
p(0,S(0)) &:=& \inf \left\{w\in \R~:~\exists\; (\phi,L)\in {\cal
A},
                     \; X_{w{\bf 1}_{1}}^{\phi,L}(T)\succeq g(S(T))
                     \right\} \;,
\end{eqnarray*}
i.e. $p(0,S(0))$ is the minimal initial capital which allows to
hedge the contingent claim $g$ by means of some admissible trading
strategy.


\section{The explicit characterization}\label{sec thmmain}

Before to state our main result, we need to define some additional
notations. We first introduce the positive polar of $K$
\begin{equation}\label{eq reec Kstar}
K^*\;:=\;\{ \xi \in \R^{1+d_c}~:~ \xi\cdot  x \geq 0, \;\; \forall x
\in K\}\;=\;\{\xi\in \R^{1+d_c}_+~:~\xi^j\le \xi^i(1+\lambda^{ij})\}\;,
\end{equation}
together with its (compact) section
  \[
    \Lambda := \{\xi \in K^*~:~\xi^1=1\}\;\subset\;(0,\infty)^{1+d_c}\;\;.
  \]
One easily checks that $\Lambda$ is not empty since it contains
the vector of $\R^{1+d_c}$ with all component equal to one. It is
a standard result in convex analysis that the partial ordering
$\succeq$ can be characterized  in terms of $\Lambda$ by
\begin{eqnarray}\label{eq cara K par Lambda}
x_1 \;\succeq\; x_2 &\mbox{if and only if}& \xi\cdot  (x_1-x_2)
\geq 0 \hspace{5mm} \mbox{for all }\; \xi\in \Lambda \;,
\end{eqnarray}
see e.g. Rockafellar (1970).

\no For $\xi \in \R^{1+d_c}$, we denote by $\lb{\xi}$ the vector
of $\R^{d_c}$ defined by $\lb{\xi}^i$ $=$ $\xi^{i+1}$ for  $i\le
d_c$. This amounts to removing the first component. With this
notations, we define
\begin{eqnarray*}
G(z)=G(z^f,z^c) &:=& \sup_{\xi \in \Lambda}\;
          \xi\cdot g\left(z^f,\mbox{diag}[\lb{\xi}]^{-1}z^c\right)
\hspace{5mm} \mbox{for  }\;z=(z^f,z^c)\mbox{ in
}\R_+^{d_f+d_c}\;,
\end{eqnarray*}
and denote by $\hat G$ the concave enveloppe of $G$ with respect
to $z^c$.

\subsection{Main result}

\no  The main result of this paper requires the additional
conditions~:

\no {$\bf (H\lambda)$ :} $\lambda^{ij}+\lambda^{ji}>0$ for all
$i,j=1,\ldots,1+d_c$, $i\neq j$.

\no {$\bf (H\sigma)$ :} For all  $1\le i\le d_f$  and   $t\le T$,
$\{\sigma(t,z)^{ij}\}_{j\le d}$ depends only on $z^f$.

\no {$\bf (Hg)$ :} $g$ is lower-semicontinuous, $\hat G$ is
continuous and has linear growth.
\\

\no The condition $(H\lambda)$ means that there is no way to avoid
transaction costs when transacting on  the $d_c$ last assets.

\begin{Theorem}\label{thmmain} Assume that $(H\lambda)$-$(H\sigma)$-$(Hg)$ hold.
Then,
    \b*
    p(0,S(0)) &=& \min \left\{w\in \R~:~\exists\; (\phi,L)\in {\cal A}^{BH},
                     \; X_{w{\bf 1}_{1}}^{\phi,L}(T)\succeq g(S(T))
                     \right\} \;,
    \e*
where
    \b*
    {\cal A}^{BH}&:=& \left\{(\phi,L)\in {\cal A}~:~L(t)=L(0) \mbox{ for all } 0\le t \le T\right\}\;.
    \e*
Moreover, there is some $\hat \Delta$ $\in$ $\R^{d_c}$ such that
    \b*
    p(0,S(0)) &=& \Esp{C(S^f(T);\hat \Delta)}
    + \sup_{\xi \in \Lambda} \lb{\xi}\cdot  \diag[\hat \Delta] S^c(0)
    \e*
where, for $\Delta \in \R^{d_c}$,
    \b*
    C(S^f(T);  \Delta):=\sup_{s^c \in (0,\infty)^{d_c}} \;
    \hat G\left(S^f(T), s^c\right)-  \Delta\cdot  s^c\;,
    \e*
 and there is an optimal hedging strategy $(\phi,L)$ $\in \Ac^{BH}$ satisfying  $L=\hat \Delta$ on $[0,T]$.
\end{Theorem}

\no The proof of the last result will be provided in the
subsequent sections.

\medskip

\no As in the papers quoted in the introduction, we obtain that the cheapest hedging strategy consists in
keeping the number of non-freely exchangeable assets, $S^c$,
constant in the portfolio. But here there is a remaining part,
namely $g(S(T))-(0,\diag[\hat \Delta] S^c(T))$, which has to be
hedged dynamically by investing in the freely exchangeable assets,
$S^f$. It is done by hedging $ C(S^f(T);\hat \Delta)$.
\\

\no From Theorem \ref{thmmain}, we can now deduce an explicit
formulation for $p(0,S(0))$.

\begin{Corollary}\label{cor formulation duale} Let the conditions
of Theorem \ref{thmmain} hold. Then,
    \b*
    p(0,S(0))
    &=&
    \min_{\Delta \in \R^{d_c}}
    \left\{
    \Esp{C(S^f(T);\Delta)}
    + \sup_{\xi \in \Lambda} \lb{\xi}\cdot  \diag[ \Delta ] S^c(0)
    \right\}
    \;.
    \e*
Moreover, if $\hat \Delta$ solves the above optimization problem,
then there is an optimal hedging strategy $(\phi,L)$  $\in
\Ac^{BH}$ which satisfies $L=\hat \Delta$ on $[0,T]$.
\end{Corollary}

\no The proof will be provided in Section \ref{sec proof finale
thmmain}.

\begin{Remark}\label{rem sur thmmain}{\rm Let the conditions
of Theorem \ref{thmmain} hold.

\no 1. In Touzi (1999), the result of Bouchard and Touzi (2000) is
generalized to the case where the initial wealth, before to be
increased by the super-replication price, is non-zero, i.e. the
following problem is considered~:
    \b*
    p(0,S(0);x)
    &:=&
    \inf \left\{w\in \R~:~\exists\; (\phi,L)\in {\cal A},
                     \; X_{x+w{\bf 1}_{1}}^{\phi,L}(T)\succeq g(S(T))
                     \right\} \;\;,
    \e*
$x \in   \R^{1+d_c}$.  Our result can be easily extended to this
case. Indeed, it suffices to observe from the wealth dynamics that
    \b*
    &X_{x+w{\bf 1}_{1}}^{\phi,L}(T)\succeq g(S(T))
                    &
    \\
    &\Longleftrightarrow&
    \\
    &X_{w{\bf 1}_{1}}^{\phi,L}(T)\succeq g(S(T))
    -\left(x^1,  \diag[S^c(0)]^{-1}\diag[S^c(T)]\lb{x}\right)
                     &
    \e*
where we recall that $\lb{x}$ is obtained from $x$ by dropping
the first component. Hence, to characterize $p(0,S(0);x)$, it
suffices to replace $g$ by
    \b*
    g(s;x)&:=& g(s) -(x^1,\diag[S^c(0)]^{-1} \diag[s^c]\lb{x})\;.
    \e*
We then deduce from Theorem \ref{thmmain} and Corollary \ref{cor
formulation duale} that, for some $\hat \Delta(x) \in \R^{d_c}$,
    \be
    p(0,S(0);x)
    &=&
    \min \left\{w\in \R~:~\exists\; (\phi,L)\in {\cal A}^{BH},
                     \; X_{x+w{\bf 1}_{1}}^{\phi,L}(T)\succeq g(S(T))
                      \right\}
    \nonumber
    \\
    &=&
    \Esp{C(S^f(T);\hat \Delta(x),x)}
    + \sup_{\xi \in \Lambda} \lb{\xi}\cdot  \diag[\hat \Delta(x)] S^c(0) \;
    \nonumber
    \\
    &=&
    \min_{\Delta \in \R^{d_c}}
    \Esp{C(S^f(T);  \Delta ,x)}
    + \sup_{\xi \in \Lambda} \lb{\xi}\cdot  \diag[  \Delta ] S^c(0)
    \;\label{eq prob min sur Delta}
    \ee
where
    \b*
      C(S^f(T);\Delta ,x):=\sup_{s^c \in (0,\infty)^{d_c}} \;
    \hat G\left(S^f(T), s^c\right)-\left(x^1, \diag[S^c(0)]^{-1}\diag[s^c]\lb{x}\right)-  \Delta \cdot  s^c\;.
    \e*

\smallskip

\no 2. The set of initial wealth which allow to hedge $g$,
    \b*
    \Gamma(g)&:=&
    \left\{x\in \R^{1+d_c}~:~\exists\; (\phi,L)\in {\cal A},\; X_{x}^{\phi,L}(T)\succeq g(S(T)) \right\}
    \;,
    \e*
can be written in
    \b*
    \Gamma(g)&=&
    \left\{x\in \R^{1+d_c}~:~p(0,S(0);x)\le 0\right\}\;.
    \e*

\smallskip

\no 3. In the limit case where $d_f=0$, we recover the result of
Bouchard and Touzi (2000) and Touzi (1999).
     }
\end{Remark}

\subsection{Example}

\no We conclude this section with a simple example. We consider a
two dimensional Black and Scholes model, i.e. $d_f=d_c=1$,
$\sigma(t,s)=\sigma \in \M^2$ with $\sigma$ invertible. In this
case, we have
    \b*
    \Lambda &=& \left\{(1,y) \in \R^2~:~
    \frac{1}{1+\lambda^{21}} \le y \le 1+ \lambda^{12}
    \right\}\;\;,\;\;\lambda^{21}+\lambda^{12}\;>\;0\;.
    \e*
We take $g$ of the form
    \b*
    g(s)&=&
    \left([s^1-K^1]^+\1_{\{s^2> K^2\}}\right)\1_1
    \;=\;
    \left([s^f-K^1]^+\1_{\{s^c> K^2\}}\right)\1_1
    \e*
with $K^1$, $K^2>0$. Then,
    \b*
    G(s)\=\left([s^f-K^1]^+\1_{\{s^c> \tilde K^2 \}}\right)\1_1
    &\And&
    \hat G(s)\= [s^f-K^1]^+ \left((s^c/\tilde K^2)\wedge 1\right)\;,
    \e*
where $ \tilde K^2=K^2/(1+\lambda^{21})$. For $\Delta \in \R$, we
have
    \b*
    C(s^f;\Delta)
    &=&
    \sup_{s^c \in (0,\infty)^{d_c}}
    \; [s^f-K^1]^+ \left((s^c/\tilde K^2)\wedge 1\right)- \Delta s^c
    \\
    &=&
    \left\{
    \begin{array}{ll}
     \left([s^f-K^1]^+ - \Delta \tilde K^2\right)
     \1_{\{0\le \Delta \tilde K^2\le[s^f-K^1]^+\}}\;\;,\;
     & \mbox{ if } \Delta\ge 0\;,
    \\
    \infty & \mbox{ otherwise.}
    \end{array}
    \right.
    \e*
Then, by Corollary \ref{cor formulation duale},
    \be
    p(0,S(0))&=& \min_{\Delta \ge 0 }
                \Esp{ C(S^f(T); \Delta) } + (1+\lambda^{12}) \Delta  S^c(0)
    \nonumber
    \\
    &=&
    \min_{\Delta \ge 0 }\;  \Esp{ \left([S^f(T)-K^1]^+ - \Delta \tilde K^2\right)
     \1_{\{0\le \Delta \tilde K^2\le[S^f(T)-K^1]^+  \}}} \nonumber
    \\
    &+& (1+\lambda^{12}) \Delta S^c(0) \nonumber
    \\
    &=&
    \min_{\Delta \ge 0 }\;  \Esp{ \left([S^f(T)-K^1- \Delta \tilde K^2]^+ \right)}
    + (1+\lambda^{12}) \Delta S^c(0)
    \;, \label{eq opti delta}
    \ee
where the expectation is convex in $\Delta$. Then, if the optimal
$\Delta$ is different from $0$, it must satisfy  the first order
condition
    \be \label{eq first order}
    -  \tilde K^2 \;\Esp{ \1_{\{S^f(T)-K^1\ge  \Delta \tilde K^2\}}}
    + (1+\lambda^{12}) S^c(0) &=& 0\;.
    \ee
We consider two different cases.

\smallskip

\no {\rm 1.} If $\Pro{S^f(T)-K^1\ge 0}$ $\le$ $(1+\lambda^{12})
S^c(0)/\tilde K^2$, then, either the only solution of \reff{eq
first order} is $0$ or \reff{eq first order} has no solution. It
follows that the optimum in \reff{eq opti delta} is achieved by
$\hat \Delta=0$. Therefore
    \b*
    p(0,S(0))
    &=& \Esp{  [S^f(T)-K^1]^+  }\;,
    \e*
 and, by the Clark-Ocone's formula,  the optimal hedging strategy $(\phi,L)$ is defined by $L=0$ and
$\phi(t)$ $=$ $ \Esp{S^f(T)\1_{\{S^f(T)\ge
K^1\}}~|~\Fc_t}/S^f(t)$.

\smallskip

\no {\rm 2.} If $\Pro{S^f(T)-K^1\ge 0}$ $>$ $(1+\lambda^{12})
S^c(0)/\tilde K^2$, then \reff{eq first order}  has a unique
solution $\hat \Delta>0$ which satisfies
    \b*
    p\;:=\;p(0,S(0))
    &=&
    \Esp{  [S^f(T)-K^1 - \hat \Delta \tilde K^2]^+  }
    + (1+\lambda^{12}) \hat \Delta S^c(0)\;.
    \e*
Observe that, in this model, $\hat \Delta$ can be computed
explicitly in terms of the inverse of the cumulated distribution
of the gaussian distribution. Let $(\phi,L)$ be defined by
    \b*
    L(t)\=\hat \Delta &\And& \phi(t)\=\Esp{S^f(T)\1_{\{S^f(T)- K^1\ge  \hat \Delta \tilde K^2\}}~|~\Fc_t}/S^f(t)
     \;\;\mbox{ on } t\le T\;.
    \e*
By the Clark-Ocone's formula, we have
    \b*
    X_{p\1_1}^{\phi,L}(T)
    &=&
    \left(
    [S^f(T)-K^1 - \hat \Delta \tilde K^2]^+\;,\; \hat \Delta S^c(T)\right)
    \;.
    \e*
For ease of notations, let us define
    \b*
    \Psi:=[S^f(T)-K^1 ]^+\1_{\{S^c(T)\ge K^2\}}\;.
    \e*
Fix $\omega \in \{\Psi\ge 0\}$. Then, $S^f(T)\ge K^1$ and
$S^c(T)\ge K^2$. If $S^f(T)-K^1 \le \hat \Delta \tilde K^2$ then
$X_{p\1_1}^{\phi,L}(T)$ $=$ $(0,\hat \Delta S^c(T))$. Recalling
the definition of $\tilde K^2$, we then obtain
    \b*
    X_{p\1_1}^{\phi,L}(T)
    \=
    (0,\hat \Delta S^c(T))
    &\succeq&
    (\hat \Delta  \tilde K^2,0)
    \;\succeq\;
    (\Psi,0)\;.
    \e*
If $S^f(T)-K^1> \hat \Delta \tilde K^2$, then
    \b*
    X_{p\1_1}^{\phi,L}(T)
    \=
    \left(S^f(T)-K^1-\hat \Delta \tilde K^2,\hat \Delta S^c(T)\right)
    &\succeq&
    (S^f(T)-K^1,0)
    \;=\;
    (\Psi,0)\;.
    \e*
On $\{\Psi=0\}$, we have $X_{p\1_1}^{\phi,L}(T) \succeq 0$ $=$
$(\Psi,0)$ since $\hat \Delta>0$.

\section{Fictitious markets}\label{sec marche fictif}

In this section, we follow the arguments of Bouchard and Touzi
(2000), i.e. we introduce an auxiliary control problem which can
be interpreted as a super-replication problem in a fictitious
market without transaction costs but were $S^c$ is replaced by a
controlled process evolving in the "bid-ask" spreads associated to
the transaction costs $\lambda$. This is obtained by introducing a
controlled process $f(Y^{(a,b)})$, see below, which evolves in
$\Lambda$. Then, the fictitious market is constructed by replacing
$S$ by $(S^f,\diag[\lb{f}(Y^{(a,b)})]S^c)$ and $g(S(T))$ by
$f(Y^{(a,b)}(T))\cdot g(S(T))$. In this paper, we shall not enter
into the detailed construction as it follows line by line the
arguments of Bouchard and Touzi (2000) up to obvious modifications
(at the level of notations). We only state the most important
results and refer to Bouchard and Touzi (2000) for the proofs, see
also the survey paper Touzi (1999).

\subsection{Parameterization of the fictitious markets}

\no We first parameterize the compact set $\Lambda$. Since $K^*$ is a polyhedral closed convex cone, we can find a family $e$ $=$ $(e_i)_{i\le n}$ in $(0,\infty)^{1+d_c}$, for some $n\ge 1$, such that, for all $\alpha \in \R^n_+$,  $ \sum_{i=1}^n \alpha^i e_i$ $=$ $0$ implies $\alpha=0$,  and $K^*$ $=$ $\{ \sum_{i=1}^n \alpha^i e_i\;,$ $\alpha \in \R^n_+\}$. Then, we define the  map $f$ from $(0,\infty)^n$ into $\Lambda$ by
    \b*
    f(y)&:=& \left(\sum_{i=1}^n y^i e_i\right)/\left(\sum_{i=1}^n y^i e_i^1\right)\;\;,\;\;y\in (0,\infty)^n\;.
    \e*
Before to go on with the definition of the fictitious markets, we list some useful properties of $f$ and $\Lambda$.

\begin{Lemma}\label{lem f} Let $(H\lambda)$ hold. Then,

\no (i) There is some $\delta>0$ such that  $0<
\xi^i+\frac{1}{\xi^i} \le \delta$ for each $\xi \in \Lambda$.

\no (ii) On  $(0,\infty)^{n}$, the rank of the Jacobian matrix $D f$ of $f$ is $d_c$.

\end{Lemma}

\proof See  Bouchard and Touzi (2000) for the easy proof.
\ep\\

\no In order to alleviate the notations, we define $\bar f$,
$\lb{F}$ and $\bar F$ as
 \b*
 \bar f^i&=& \lb{f}^{i-d_f} \1_{\{i> d_f\}} + \1_{\{i\le d_f\}}\;\;,\;\;
 \lb{F}\;=\;\diag[\lb{f}] \;\;\;\And\; \bar F\;=\;\diag[\bar f]
 \;.
 \e*
The map $\bar f$ coincides with $f$ on its last $d_c$ components,
while the first $d_f$ ones are set to one. Here, $\bar f$, $\bar
F$ and $\lb{F}$ take values in $\R^d$, $\M^d$ and $\M^{d_c}$
respectively.\\

\no Given some arbitrary parameter $\mu
>0$, we then define for all $(y_{0},z_0) \in (0,\infty)^{n} \x \R_+^d$ the
continuous function $\alpha ^{y_{0},s_0}$ on $[0,T]\times
\R_+^d\x (0,\infty)^{n}\times \M^{n,d}\times \R^{n}$ as
\begin{equation}
\alpha^{y_{0},s_0}(t,s,y,a,b):=\left\{
\begin{array}{ll}
A(t,s,y,a,b) & \mbox{ \quad if } \sum_{i=1}^{d}\sum_{j=1}^{n}
\left( |  s^{i}-z^{i}_0 |+|\ln
\frac{y^{j}}{y_{0}^{j}}|\right) <\mu \mbox{
\quad } \\
\mbox{constant}&\mbox{otherwise,}
\end{array}
\right.  \label{alpha}
\end{equation}
where
\begin{eqnarray*}
A(t,s,y,a,b)
&=\; {\sigma}(t,s)^{-1}  \bar F(y)^{-1}&
\left\{D \bar f(y)\diag[y]b \right. \\
&& + \frac{1}{2}\mbox{Vect}\left[\mbox{Tr}\left( D^2 \bar
f^{i}(y)\diag[y]aa'
                                         \diag[y]\right)\right]_{i\le d} \\
&& + \left.\mbox{Vect}\left[ \left(D \bar
f(y)\diag[y]a\sigma(t,s)'\right)_{ii }\right]_{i\le d} \right\}\;.
\end{eqnarray*}
Let $\Dc$ be the set of all bounded progressively measurable
processes $(a,b)=\{(a(t),b(t)),$ $0\leq t\leq T\}$ where $a$ and
$b$ are valued respectively in $\M^{n,d}$ and $\R^n$. For all
$(t,y,z)$ in $[0,T]\x (0,\infty )^{n+d}$ and $(a,b)$ in ${\cal
D},$ we introduce the controlled process $Y_{t,y,z}^{(a,b)}$
defined on $[t,T]$   as the solution of the stochastic
differential equation
    \be
    dY(r) &=&\mbox{diag}[Y(r)]\left[ \left(
    b(r)+a(r)\alpha ^{y,s}(r,S_{t,s}(r),Y(r),a(r),b(r))\right)
    dt+a(r)dW(r)\right] \nonumber\\
    Y(t) &=& y\;,\label{Y_ab}
    \ee
where $S_{t,s}$ is the solution of \reff{def_s} with the condition
$S_{t,s}(t)=s$ and $s=\bar F(y)^{-1}z$. It follows from our
assumption on $\sigma$  that
$\alpha^{y_{0},s_0}(t,s,y,a,b)$ is a random Lipschitz function of
$y$, so that the process $Y_{t,y,z}^{(a,b)}$ is well defined on
$[t,T].$ For each $(a,b)$ in ${\cal D},$ we define the process
$Z_{t,y,z}^{(a,b)}$ by
\begin{equation}
Z_{t,y,z}^{(a,b)}=\bar F(Y_{t,y,z}^{(a,b)})S_{t,s} \;\;\mbox{ with }\; s=\bar F(y)^{-1}z\;.
\label{Z_ab}
\end{equation}
Observe that $(Z_{t,y,z}^{(a,b)})^f$ $=$ $S_{t,s}^f$.

\subsection{Super-replication in the fictitious markets}

Let $\phi$ be a progressively measurable process valued in
$\R^{d}$  satisfying
    \begin{equation}
    \sum_{i=1}^{d} \int_{0}^{T}|\phi ^{i}(t)|^{2}d\langle Z_{t,y,z}^{(a,b),i}(t)\rangle <\infty\;.
    \label{integrabilite_phi}
    \end{equation}
Then, given $x\geq 0,$ we introduce the process $X_{t,x,y,z}^{(a,b)^{\phi}}$
defined by
    \begin{equation}
    X_{t,x,y,z}^{(a,b)^{\phi}}(r)=x+\int_{t}^{r}\phi (s)\cdot
    dZ_{t,y,z}^{(a,b)}(s) \label{W_ab}
    \end{equation}
and we denote by ${\cal B}^{(a,b)}(t,x,y,z)$ the set of all such processes $\phi$
satisfying the additional condition
    \begin{equation}\label{phiadm}
    X_{t,x,y,z}^{(a,b)^{\phi}}(r)
    \;\geq\;
    -c-\delta\cdot  Z_{t,y,z}^{(a,b)}(r) \hspace{2mm}
    , \hspace{2mm} t\leq r\leq  T\;,\;\mbox{ for some } (c,\delta) \in \R^{1+d}\;.
    \end{equation}
We finally define the auxiliary stochastic control problems
    \begin{eqnarray}
    u^{(a,b)}(t,y,z)
    &:=&
    \inf \left\{ x\in \R~:~\exists \phi \in {\cal B}^{(a,b)}(t,x,y,z) \;,\right.
    \nonumber\\
    & & \hspace{7mm} \left.
    X_{t,x,y,z}^{(a,b)^{\phi}}(T)\geq f\left(Y_{t,y,z}^{(a,b)}(T)\right)\cdot  g\left( S_{t,s}(T)\right)
    \right\} \;, \label{u_ab}
    \end{eqnarray}
with $s$ $=$ $\bar F(y)^{-1}z$, and
\begin{equation}
u(t,y,z):=\sup_{(a,b)\in {\cal D}}u^{(a,b)}(t,y,z) \;.\label{u}
\end{equation}

\no The value function $u^{(a,b)}(t,y,z)$ coincides with the
super-replication price of the modified claim $f\left(
Y_{t,y,z}^{(a,b)}(T)\right)\cdot  g\left( S_{t,s}(T)\right)$ in
the market formed by the assets $Z_{t,y,z}^{(a,b)}$ without
transaction costs. The function $u(t,y,z)$ is the upper-bound of
these prices over all the ``controlled" fictitious markets. We
refer to Bouchard and Touzi (2000) for a more detailed discussion.

\subsection{Viscosity properties of $u_*$}

\no We can now provide a first lower bound for $p(0,S(0))$ which
is similar to the one provided by Bouchard and Touzi (2000) in
the case $d_f=0$.

\no For $(t,y,z) \in [0,T]\x (0,\infty)^n\x\R_+^d$, we define the lower semicontinuous function $u_*$ by
    \b*
    u_*(t,y,z)&:=& \liminf_{{\hspace{5mm} }^{(t',y',z')\to (t,y,z)}_{   (t',y',z')\in [0,T)\x (0,\infty)^{n+d}} } u(t',y',z')\;.
    \e*
Contrary to Bouchard and Touzi (2000), we need to extend the
definition of $u_*$ to the whole subspace $[0,T]\x
(0,\infty)^n\x\R_+^d$ (in opposition to $[0,T]\x
(0,\infty)^{n+d}$). Although, we are only interested by $u_*$ on
$[0,T)\x (0,\infty)^{n+d}$, since $S(0) \in (0,\infty)^{d}$, this
extension will be useful to apply the comparison theorem of
Proposition \ref{prop theo comp} below.

\begin{Theorem}\label{thm BT} Let $(H\lambda)$ and $(Hg)$ hold. Then $u_*$ satisfies~:

\no (i) For all $y$ $\in$ $(0,\infty)^{n}$, we have $p(0,S(0)) \ge
u_*(0,y,\bar F(y)S(0))$.

\no (ii)  $u_*$ is independent of its variable $y$.

\no (iii) $u_*$ is a viscosity supersolution on
$[0,T)\x(0,\infty)^{n}\x\R_+^d$ of
    \b*
    \inf_{a \in \M^{n,d}}
    -\Hc^a \vp \ge 0 \;,
    \e*
where, for a smooth function $\vp$,
    \b*
    \Hc^a\vp&:=& \frac{\partial \vp}{\partial t} + \frac{1}{2} \Tr\left[\Gamma^{a'} D^2_{zz} \vp \Gamma^a\right]
    \e*
and
    \b*
    \Gamma^{a}(t,y,z)
    &:=&
    \diag[z]
    \left( \sigma(t,\bar F(y)^{-1}z)+\bar F(y)^{-1}D  \bar f(y)\diag[y]a \right)
    \;.
    \e*
\no (iv) For all $(y,z) \in (0,\infty)^{n}\x \R_+^d$
    \b*
    u_*(T,y,z)&\ge& G(z)\;.
    \e*
\end{Theorem}

\no This result is obtained by following line by line the
arguments of Sections 6, 7 and 8 in Bouchard and Touzi (2000),
see also Touzi (1999). Since its proof is rather long, we omit it.
\\

\no In Bouchard and Touzi (2000), the above characterization was
sufficient to solve the super-replication problem. Indeed, in the
case where $d_f=0$, one can show that $u_*$ is concave with
respect to $z$ and  non-increasing in $t$. This turns out to be
sufficient to show that it corresponds to the price of the
cheapest buy-and-hold super-hedging strategy in the original
market. In our context, where $d_f\ge 1$, we can only show that
$u_*$ is concave with respect to $z^c$ and there is no reason why
it should be concave in $z$ (in particular   if $g(s^f,s^c)$
depends only on $s^f$). We therefore have to work a little more.

\no  As a first step, we rewrite the above PDE in a more tractable
way. For all $(t,z) \in[0,T] \times \R_+^d$ and $\mu \in
\M^{d_c,d}$, we  define
    \b*
    \sigma^\mu(t,z)&:=& \diag[z]\left[\sigma(t,z)^{ij}\1_{i\le d_f} +\mu^{ij}1_{i> d_f}
    \right]_{1\le i,j\le d}
    \e*
where, for real numbers $(a^{ij})$, $\left[a^{ij}\right]_{1\le
i,j\le d}$ denotes the square $d$-dimensional matrix $M$ defined
by $M^{ij}=a^{ij}$.

\no Since $u_*$ does not depend on its $y$ variable, from now on,
we shall omit it if not required by the context.

\begin{Corollary}\label{cor BT} Let $(H\lambda)$, $(H\sigma)$ and $(Hg)$ hold.
Then,

\no (i) $u_*$ is a viscosity supersolution on
$[0,T)\x\R_+^d$ of
    \be\label{eq viscosite inf sur U}
    \inf_{\mu \in \M^{d_c,d}}
    -\Gc^\mu\vp \ge 0 \;,
    \ee
where, for a smooth function $\vp$ and $\mu \in \M^{d_c,d} $,
    \b*
    \Gc^\mu\vp(t,z)&=&
    \frac{\partial \vp}{\partial t}
    +
    \frac{1}{2} \mbox{Tr}\left[\sigma^\mu(t,z)'  D^2_{zz}\vp(t,z) \sigma^\mu(t,z)\right]
    \;.
    \e*
\no (ii) For each $(t,z^f)\in [0,T)\x(0,\infty)^{d_f}$, the map
$z^c\in (0,\infty)^{d_c}\mapsto u_*(t,z^f,z^c)$ is concave.

\no (iii) For all  $z\in \R_+^d$
    \be\label{eq viscosite inf sur U bord}
    u_*(T,z)&\ge& \hat G(z)\;,
    \ee
where we recall that $\hat G$ is the concave envelope of $G$ with
respect to its last $d_c$ variables.
\end{Corollary}

\proof (i). Recall from Lemma \ref{lem f} that the rank of $D
f(y)$ is $d_c$ whenever $y \in (0,\infty)^{n}$. Since $\bar f^i$ $=$ $1$
for $i\le d_f$, we deduce from $(H\sigma)$ that, for each $\mu \in \M^{d_c,d}$, we
can find some $a$ $\in$ $\M^{n,d}$ such that $\Gamma^a(t,y,z)$
$=$ $\sigma^\mu(t,z)$. Then, the first result follows from
Theorem \ref{thm BT}.

\no (ii). For $\vp$ satisfying \reff{eq viscosite inf sur U} we
must have, on $[0,T)\x(0,\infty)^d$, $-Tr[\mu' D^2_{z^cz^c} \vp \mu]\ge 0$ for all $\mu \in
\M^{d_c,d}$ since otherwise we would get a contradiction of
 \reff{eq viscosite inf sur U} by considering $\delta \mu$ and
sending $\delta$ to infinity. Then, the concavity property follows
from the same argument as in Lemma 8.1 of Bouchard and Touzi (2000).

\no (iii). In view of the boundary condition of   Theorem \ref{thm BT}, is suffices to show
 that $u_*(T,z)$ is concave with respect to $z^c$. To see this,
 fix $(z^f,z^c_1,z^c_2)$ $\in  \R_+^{d_f+2d_c}$ and
 observe that, by (ii) and definition of $u_*$,
    \b*
    u_*(T, z^f,(z^c_1+z^c_2)/2)
    &=&
    \liminf_{^{(\tilde t, \tilde z^f,\tilde z^c_1,\tilde z^c_2)\to (T, z^f,z^c_1,z^c_2)}_{(\tilde t, \tilde z^f,\tilde z^c_1,\tilde z^c_2) \in [0,T)\x          (0,\infty)^{d_f+2d_c}}}
    u_*(\tilde t, \tilde z^f,(\tilde z^c_1+\tilde z^c_2)/2)
    \\
    &\ge&
    \liminf_{^{(\tilde t, \tilde z^f,\tilde z^c_1,\tilde z^c_2)\to (T, z^f,z^c_1,z^c_2)}_{(\tilde t, \tilde z^f,\tilde z^c_1,\tilde z^c_2) \in [0,T)\x          (0,\infty)^{d_f+2d_c}}}
    \frac{1}{2}
    \left(
    u_*(\tilde t, \tilde z^f, \tilde z^c_1  ) + u_*(\tilde t, \tilde z^f, \tilde z^c_2  )
    \right)
    \\
    &\ge&
    \frac{1}{2}
    \left(u_*(T,   z^f,   z^c_1  ) + u_*( T,   z^f,   z^c_2 )\right)\;.
    \e*
 \ep


\section{A tractable lower bound for $u_*$}\label{sec pbm auxiliaire}

\no Let us introduce some additional notations.
Given $\kappa\ge 0$, we define
$U_\kappa$ as the set of all elements $M$ of $\M^{d_c,d}$
such that $|M|\le \kappa$. We then denote by $\Uc_\kappa$  the collection of all
$U_\kappa$-valued predictable processes. To each $\mu \in
\bigcup_{\kappa\ge 0} \Uc_{\kappa}$, we   associate the
controlled process $Z^\mu_{t,z}$ defined as the solution of
    \be\label{eq def Zmu}
    Z(s) &=& z+ \int_t^s\sigma^{\mu(r)}(r,Z(r)) dW(r)\;\;\;\; t\le s \le
    T\;.
    \ee
Observe that, under $(H\sigma)$,
    \b*
    (Z^\mu_{t,z})^f&=& S_{t,z}^f \;.
    \e*

\no The aim of this Section is to prove the following result.

\begin{Proposition}\label{prop u ge w} Let $(H\lambda)$, $(H\sigma)$ and  $(Hg)$ hold. Then,
for all $z \in (0,\infty)^{d}$,
    \b*
        u_*(0,z)&\ge&   \sup_{\mu \in \Uc} \Esp{\hat    G(Z_{0,z}^\mu(T))}\;,
    \e*
where $\Uc$ denotes the set of all $\M^{d_c,d}$-valued square
integrable predictable processes $\mu$ such that $Z_{0,z}^\mu$ is
a  martingale for all $z$ $\in$ $(0,\infty)^{d}$.
\end{Proposition}

\no Using an approximation argument combined with the martingale
property of the $Z^\mu$'s, the concavity of $u_*$ with respect to
$z^c$ and  assumption $(H\sigma)$,  this will allow us to show
that,  for all $z=(z^f,z^c)\in (0,\infty)^{d_f} \x
(0,\infty)^{d_c}$ and  $\Delta\in \partial_{z^c}u_*(0,z^f,z^c)$,
we have
    \b*
    u_*(0,z^f,z^c)
    &\ge&
    \Esp{
    \sup_{\tilde z^c \in (0,\infty)^{d_c}}
    \left\{\hat  G(S^f_{0,z}(T),\tilde z^c)-\Delta\cdot  \tilde z^c\right\}
    }
    +\Delta\cdot  z^c\;,
    \e*
where,
$\partial_{z^c}u_*(0,z^f,z^c)$ is the subgradient of the mapping
$z^c\mapsto u_*(0,z^f,z^c)$, see Corollary \ref{cor borne finale
pour u} below. This last lower bound for $u_*$ will turn out to be
enough to conclude the proof, see Section \ref{sec proof finale thmmain}. \\

\no To this purpose, we shall first consider the auxiliary control problems
    \be
    v_\kappa(t,z)&:=& \sup_{\mu \in
    \Uc_{\kappa}}\Esp{\hat G\left(Z^\mu_{t,z}(T)\right)}\;\;\;\;(t,z,\kappa) \in
    [0,T]\x (0,\infty)^{d} \x (0,\infty)\;,
    \ee
 and show that $u_*\ge \sup_{\kappa>0} v^*_\kappa$, where $v^*_\kappa$ is the upper-semicontinuous function defined on $[0,T]\x \R^d_+$ by
    \b*
    v_\kappa^*(t, z)&:=& \limsup_{{\hspace{5mm} }^{(t' ,z')\to (t ,z)}_{   (t' ,z')\in [0,T)\x (0,\infty)^{d}} } v_\kappa(t' ,z')\;.
    \e*
This will be done by means of a comparison argument on the
 PDE defined by \reff{eq viscosite inf sur U}-\reff{eq viscosite inf sur U bord} with $U_\kappa$
 substituted to $\M^{d_c,d}$.

\no In the next subsection, we show that $v^*_\kappa$ is a
  viscosity subsolution of \reff{eq viscosite inf
sur U}-\reff{eq viscosite inf sur U bord} with $U_\kappa$
substituted to $\M^{d_c,d}$. Then, we provide the comparison
theorem. We conclude  the proof of Proposition \ref{prop u ge w}
in the last subsection.

\subsection{Viscosity properties of $v_\kappa^*$}

\no We start with the subsolution property in the interior of the
domain. The proof is rather standard now but, as it is short, we
provide it for completeness.

\begin{Lemma}\label{lem vkappa soussol} For each $\kappa>0$, $v^*_\kappa$ is a viscosity subsolution
on $[0,T)\times \R_+^{d}$ of
    \b*
    \inf_{\mu \in U_\kappa } -\Gc^\mu \vp &\le& 0 \;.
    \e*
\end{Lemma}

\proof Let $\varphi$ $\in$ $C^{2}([0,T] \times \R^{d})$ and
$(t_0,z_0)$ be a strict global maximizer of $v_\kappa^*-\varphi$ on $[0,T)\x \R_+^{d}$ such that
$(v_\kappa^*-\varphi)(t_0,z_0)=0$. We  assume that
        \begin{eqnarray}\label{eq soussol int contra}
        \inf_{\mu \in U_\kappa } -\Gc^\mu \vp(t_0,z_0) & >& 0 \;,
        \end{eqnarray}
and work towards a contradiction. If \reff{eq soussol int contra}
holds, then it follows from our continuity assumptions on $\sigma$
that  there exists some $t_0<\eta<T-t_0$  such that
    \begin{eqnarray}
    \inf_{\mu \in U_\kappa } -\Gc^\mu \vp(t,z) &\ge& 0 \;\;\;\pourtout (t,z) \in B_0:=B((t_0,z_0),\eta) \;.
    \label{visco_subsol_contra}
    \end{eqnarray}
Recall that $B((t_0,z_0),\eta)$ is the open ball of radius $\eta$
centered on $(t_0,z_0)$, see the notations section.  Let
$(t_n,z_n)_{n\geq 0}$ be a sequence in $ B_0\cap ([0,T)\x(0,\infty)^d)$ such that
    \begin{eqnarray*}
    (t_n,z_n) \;\longrightarrow\; (t_0,z_0) &\mbox{and}&
    v_\kappa(t_n,z_n) \;\longrightarrow\; v_\kappa^*(t_0,z_0) \;\;
    \end{eqnarray*}
 and notice that
    \begin{eqnarray} \label{eu-vpto0}
    v_\kappa(t_n,z_n)-\vp(t_n,z_n) \;\lra\; 0    \;.
    \end{eqnarray}
Next, define the stopping times
    \begin{eqnarray*}
    \theta_n^\mu
    &:=&
    T\wedge\inf\left\{ s>t_n~:~(s,Z^\mu_n(s)) \not\in B_0 \right\}
    \end{eqnarray*}
where $\mu$ is any element of $\Uc_\kappa$ and $Z^\mu_n$ $:=$
$Z^\mu_{t_n,z_n}$. Let $\partial_p B_0=[t_0,t_0+\eta]\x
\partial  B(z_0,\eta) \cup \{t_0+\eta\}\x  B(z_0,\eta) $
 denote the parabolic boundary of $B_0$ and observe that
    \b*
    0\;>\;-\zeta&:=& \sup_{(t,z) \in \partial_p B_0 \cap ([0,T]\x\R^d_+) } (v_\kappa^*-\varphi)(t,z)
    \e*
since $(t_0,z_0)$ is a strict maximizer of $v_\kappa^*-\varphi$.
Then, for a fixed $\mu \in \Uc_\kappa$, we deduce from It\^o's Lemma
and   \reff{visco_subsol_contra}  that
    \begin{eqnarray*}
    \vp(t_n,z_n)
    &\ge& \Esp{\vp(\theta^\mu_n,Z_n^\mu(\theta^\mu_n))}
    \ge \Esp{v_\kappa^*\left(\theta^\mu_n,Z_n^\mu(\theta^\mu_n)\right)+\zeta}
    \ge \zeta + \Esp{\hat G(Z_n^\mu(T))}\;,
    \end{eqnarray*}
where we used the fact that $\vp\ge v^*_\kappa\ge v_\kappa$ and
    \b*
    v_\kappa\left(\theta^\mu_n,Z_n^\mu(\theta^\mu_n)\right)
    &\ge&
    \Esp{\hat G(Z_n^\mu(T))~|~\Fc_{\theta^\mu_n}}\;.
    \e*
By arbitrariness of $\mu \in \Uc_\kappa$, it follows from the
previous inequality  that
    \begin{eqnarray*}
    \vp(t_n,z_n)
    &\ge&  \zeta + v_\kappa(t_n,z_n)\;.
    \end{eqnarray*}
In view of \reff{eu-vpto0}, this leads to a contradiction since
$\zeta>0$. \ep\\

\no We now turn to the boundary condition.

\begin{Lemma}\label{lem vkappa soussol cond bord} Under $(Hg)$, for each $\kappa>0$ and $z \in \R_+^{d}$,
$v^*_\kappa(T,z)  \le \hat G(z)$.
\end{Lemma}

\proof For ease of notations, we write $\bar v_\kappa(z)$ for
$v^*_\kappa(T,z)$. Let $f$ be in $C^2(\R^d)$ and $z_0 \in
\R^d_+$ be such that
\begin{eqnarray*}
0&=&(\bar v_\kappa-f)(z_0)=\max_{\R_+^{d}}(\bar
v_\kappa-f)\;.
\end{eqnarray*}
We assume that
    \begin{eqnarray} \label{eqf>G}
    \bar v_\kappa(z_0)-\hat G(z_0)\;=\; f(z_0)-\hat G(z_0) &>& 0 \;,
    \end{eqnarray}
and work towards a contradiction to the definition of $v_\kappa$.\\

\no {{1. }} Define  on $[0,T]$ $\x$ $\R^{d}$
    \begin{eqnarray*}
    \vp(t,z)&:=& f(z)+ c|z-z_0|^2+(T-t)^{\frac{1}{2}}
    \end{eqnarray*}
and notice that for all $z$ $\in$ $\R^{d}$
    \begin{eqnarray}\label{eq dt to infty}
    \frac{\partial \vp}{\partial t}(t,z)  \;\lra\; -\infty &\mbox{as}&   t \lra T \;.
    \end{eqnarray}
Since $U_\kappa$ is compact,   there is some $\eta,\;\tilde \eta>0$ such that
    \begin{eqnarray}\label{eqLvp>0}
    \inf_{\mu \in U_\kappa} -\Gc^\mu\vp(t,z) &>& 0\quad
    \mbox{ for all } (t,z) \in B_0 := [T-\tilde \eta,T)\x \bar B(z_0,\eta)\;.
    \end{eqnarray}
Since  $\vp(T,z_0)$ $=$ $f(z_0)$ and $\hat G$ is continuous, see $(Hg)$, it
follows from  \reff{eqf>G} that we can choose $\eta$, $\tilde \eta$ such that, for
some $\eps>0$, we also have
    \begin{eqnarray}\label{eqvp>Gamma+eps}
    \vp(T,z) -\hat G(z)&>&\eps \;\;\;\;\pourtout\;z\in  \bar B(z_0,\eta)\cap  \R^d_+\;,
    \end{eqnarray}
and, by upper-semicontinuity of $ v_\kappa^*-f$,
    \be\label{eq v le bar v + alpha}
    v_\kappa^*(t,z)\;\le\; f(z) + \alpha &\mbox{ for some } \alpha>0 \mbox{ on } B_0 \cap ([0,T]\x \R^d_+)\;.
    \ee
By definition of $\vp$, we also have
    \be\label{eq vp ge f + c eta}
    \vp(t,z)&\ge& f(z) + c\eta^2 \;\;\mbox{ on } [T-\tilde \eta,T]\x \partial  \bar B(z_0,\eta) \;,
    \ee
where, by possibly taking a smaller $\tilde \eta$, we can choose $c$ large enough so that
    \be\label{eq c gd}
    c\eta^2 &\ge & \alpha +\eps \;
    \ee
and \reff{eqLvp>0}-\reff{eqvp>Gamma+eps}-\reff{eq v le bar v +
alpha} still holds, see \reff{eq dt to infty}.

\no {{2. }}  Let $(s_n,\xi_n)$ be a sequence in $[T-\tilde
\eta/2,T)\x$ $  (\bar B(z_0,\eta)\cap (0,\infty)^d) $ $\subset$ $B_0$ satisfying
    \begin{eqnarray*}
    (s_n,\xi_n) \longrightarrow (T,z_0)\;,\;\; s_n<T &\mbox{and}&
    v_\kappa^*(s_n,\xi_n) \longrightarrow \bar v_\kappa(z_0)\;.
    \end{eqnarray*}
Let $(t_n,z_n)$ be a maximizer of $(v_\kappa^*-\vp)$ on
$[s_n,T]\x$ $(\bar B(z_0,\eta)\cap [0,\infty)^d)$ $\subset$ $B_0$.
 For all $n$, let $(t_n^k,z_n^k)_k$ be a subsequence in $[s_n,T)\x (\bar B(z_0,\eta)\cap (0,\infty)^d)$
satisfying
    \begin{eqnarray*}
    (t_n^k,z_n^k) \longrightarrow (t_n,z_n) &\mbox{and}&
    v_\kappa(t_n^k,z_n^k) \longrightarrow v_\kappa^*(t_n,z_n)\;.
    \end{eqnarray*}
We shall prove later that
    \begin{eqnarray}\label{con_u^*Gb}
     (t_n,z_n) \longrightarrow (T,z_0) &\mbox{and}&  v_\kappa^*(t_n,z_n) \longrightarrow \bar v_\kappa(z_0)
    \end{eqnarray}
and that there exists a subsequence of $(t_n^k,z_n^k)_{k,n}$,
relabelled $(t'_n,z'_n)$, satisfying
    \begin{eqnarray}
     \hspace{-5mm} (t'_n,z'_n) \rightarrow (T,z_0) &\mbox{and}&  v_\kappa(t'_n,z'_n) \rightarrow \bar v_\kappa(z_0)
    \;\;  \mbox{, where }   t'_n<T \mbox{ for all $n$ }  \;.
    \label{eq_t'<T}
    \end{eqnarray}

\no {{3. }} For all $n$, we define the stopping
times
    \begin{eqnarray*}
    \theta_n^\mu
    &:=&
    T\wedge\inf\left\{ s>t'_n~:~(s,Z^\mu_n(s)) \in \partial B_0\right\}  \;\;,\;\;\mu \in U_\kappa\;,
    \end{eqnarray*}
where $Z^\mu_n$ $:=$ $Z^\mu_{t'_n,z'_n}$, together with $\Jc^\mu_n:=\{\theta_n^\mu<T\}$.  \\

\no {{4. }} We can now prove the required contradiction. For fixed
$\mu \in \Uc_\kappa$, we deduce from It\^o's Lemma and
\reff{eqLvp>0}  that
    \begin{eqnarray*}
    \vp(t'_n,z'_n)
    &\ge& \Esp{\vp(\theta^\mu_n,Z_n^\mu(\theta^\mu_n))}\;.
    \end{eqnarray*}
Recalling \reff{eqvp>Gamma+eps}, \reff{eq vp ge f + c eta}, \reff{eq v le bar v + alpha},  we obtain
    \begin{eqnarray*}
    \vp(t'_n,z'_n)
    &\ge&
    \Esp{\left(\hat G(Z_n^\mu(T))+\eps\right)\1_{(\Jc^\mu_n)^c}
    + \left(  c\eta^2-\alpha+  v^*_\kappa(\theta^\mu_n,Z_n^\mu(\theta^\mu_n))\right)\1_{\Jc^\mu_n} }\;.
    \end{eqnarray*}
Since $v_\kappa(T,\cdot)=\hat G(\cdot)$, we deduce from the previous inequality and \reff{eq c gd} that
    \begin{eqnarray*}
    \vp(t'_n,z'_n)
    &\ge& \Esp{ \eps+  v_\kappa(\theta^\mu_n,Z_n^\mu(\theta^\mu_n)) }
    \;\ge\;
    \eps  + \Esp{ \hat G(Z_n^\mu(T)) }\;.
    \end{eqnarray*}
By \reff{eq_t'<T} and definition of $\vp$, we can also choose $n$ such that that
$v_\kappa(t'_n,z'_n)\ge \vp(t'_n,z'_n)-\eps/2$, so that
    \begin{eqnarray*}
    v_\kappa(t'_n,z'_n)  &\ge& \eps/2+ \Esp{ \hat G(Z_n^\mu(T)) }
    \end{eqnarray*}
which, by arbitrariness of $\mu \in \Uc_\kappa$, contradicts the definition of $v_\kappa(t'_n,z'_n)$.\\

\no {{5. }} It remains to prove \reff{con_u^*Gb} and
\reff{eq_t'<T}.
 Clearly, $t_n \rightarrow T$. Let $\hat z \in  \bar B(z_0,\eta)\cap [0,\infty)^d $ be such that
$z_n \rightarrow \hat z$, along some subsequence. Then, by
definition of $f$ and $z_0$, we have
\begin{eqnarray*}
0 & \geq & (\bar v_\kappa-f)(\hat z)-(\bar v_\kappa-f)(z_0) \\
 & \geq& \limsup_{n \to \infty} (v^*_\kappa-\vp)(t_n,z_n) +c|\hat z-z_0|^2
 -(v^*_\kappa-\vp)(s_n,\xi_n)\\
& \geq &  c |\hat z-z_0|^2 \geq 0 \;,
\end{eqnarray*}
where the third inequality is obtained by definition of
$(t_n,z_n)$. Then, $\hat z=z_0$ and, by continuity of $\vp$ and
definition of $(s_n,\xi_n)$, $v^*_\kappa(t_n,z_n) \rightarrow \bar
v_\kappa(z_0)$. This also proves that
\begin{eqnarray}
 \lim_n\lim_k(t_n^k,z_n^k) =(T,z_0) &\mbox{and}&  \lim_n\lim_k v_\kappa(t_n^k,z_n^k)= \bar v_\kappa(z_0)
    \label{eq con_v kappa} \;.
\end{eqnarray}
Now assume that  card$\{(n,k)\in \N\x \N~:~t_n^k=T\}=\infty$.
Since $v_\kappa(T,\cdot)=\hat G(\cdot)$ and $\hat G$ is
continuous, there exists a subsequence, relabelled
$(t_n^k,z_n^k)$, such that
$$
\limsup_n \limsup_k v_\kappa(t_n^k,z_n^k) \leq   \hat G(z_0) \;.
$$
 Since by
assumption $\hat G(z_0)<f(z_0)=\bar v_\kappa(z_0)$, this leads to
a contradiction with \reff{eq con_v kappa}. Hence, card$\{(n,k)\in
\N\x \N ~:~t_n^k=T\}< \infty$, and, using a diagonalization
argument, we can construct a subsequence $(t'_n,z'_n)_n$ of
$(t_n^k,z_n^k)_{n,k}$ satisfying \reff{eq_t'<T}. \ep
\subsection{The comparison theorem}

\begin{Proposition}\label{prop theo comp}

Let $V$ be an upper semicontinuous viscosity subsolution   and
$U$ be a lower semicontinuous viscosity supersolution on
$[0,T)\x\R_+^d$ of
 \begin{eqnarray}
 \nonumber
    \inf_{\mu \in U_\kappa} -\Gc^\mu\varphi&=&0\;.
 \end{eqnarray}
Assume that $V$ and $U$ satisfy the linear growth condition
 \begin{eqnarray}
 \nonumber
   |V(t,z)|+|U(t,z)|\le
   K\,\left(1+|z|\right)&,\;(t,z)\in[0,T)\x\R_+^d\;,&K
   > 0\;.
 \end{eqnarray}
Then,
 \begin{eqnarray}
 \nonumber
    V(T,.)\;\le\; U(T,.)&\mbox{implies}&
    V\;\le\;U~~\mbox{on}~[0,T]\x\R_+^d\;.
 \end{eqnarray}
\end{Proposition}

\proof 1. Let $\lambda$ be some positive parameter and consider
the functions
 \b*
 u(t,z)\;:=\; e^{\lambda\,t}U(t,z)&\mbox{and}&
 v(t,z)\;:=\;e^{\lambda \,t}V(t,z)\.
 \e*
 It is easy to verify that the functions $u$ and $v$ are,
 respectively, a lower semicontinuous viscosity supersolution and
 an upper semicontinuous viscosity subsolution on
 $[0,T)\x\R_+^d$ of
  \begin{eqnarray}
  \label{eq proof thm comparison 1}
  \lambda \varphi -\frac{\partial \varphi}{\partial t} - \sup_{\mu \in
  U_\kappa}Tr\left[\sigma^{\mu\,\prime} D^2_{zz}\varphi
  \sigma^\mu\right]\;=\;0\;.
  \end{eqnarray}
  Moreover $u$ and $v$ satisfy
  \begin{eqnarray}
  \nonumber
  u(T,z) \;\ge\;v(T,z) &\mbox{for all}~~z\,\in\, \R_+^d\;,
  \end{eqnarray}
  as well as the linear growth condition
  \begin{eqnarray}
  \label{eq proof thm comparison 2}
    |v(t,z)|+|u(t,z)|\le
    A\,\left(1+|z|\right)&,\;(t,z)\in[0,T)\x\R_+^d\;,&A > 0\;.
  \end{eqnarray}
  Through the following steps of the proof we are going to show
  that $u$ $\ge$ $v$ on the entire domain $[0,T]\x\R_+^d$,
  which is equivalent to $U$ $\ge$ $V$ on $[0,T]\x\R_+^d$.

 \vspace{2mm}

 \no We argue by contradiction, and assume that for some $(t_0,z_0)$ in $[0,T]\x\R_+^d$
 \b*
 0\;<\; \delta\;:=v(t_0,z_0)-u(t_0,z_0)\;.
 \e*

\vspace{2mm}
 \no 2. Following Barles et al. (2003), we introduce the following functions. For some positive parameter $\alpha$, we set
  \b*
    \phi_{\alpha}(z,z^\prime)\;=\;\left[1+|z|^2\right]\left[\eps + \alpha
    |z^\prime|^2\right]
    &\mbox{and}&
    \Phi_{\alpha}(t,z,z^\prime)\;=\; e^{L (T-t)}\phi_{\alpha}(z+z^\prime,z-z^\prime)\;.
  \e*
 Here, $L$ and $\eps$ are positive constants to be chosen later and we don't write the
 dependence of $\phi_{\alpha}$, $\Phi_{\alpha}$ and $\Psi_{\alpha}$ with respect to them.

 \no By the linear growth condition \reff{eq proof thm comparison 2}, the upper semicontinuous function
 $\Psi_{\alpha}$ defined by
 \b*
 \Psi_{\alpha}(t,z,z^\prime)\;:=\;
 v(t,z)-u(t,z^\prime) - \Phi_{\alpha}(t,z,z^\prime)
 \e*
 is such that for all $(t,z,z^\prime)$ in $[0,T]\x\R_+^{2d}$
 \b*
 \Psi_{\alpha}(t,z,z^\prime)&\le&
  A\left(1+ |z|+ |z^\prime|\right) -
  \min\left\{\eps,\alpha\right\}
  \left(|z-z^\prime|^2+|z+ z^\prime|^2+1\right)\\&\le&
  A\left(1+ |z|+ |z^\prime|\right) - \min\left\{\eps,\alpha\right\}
  \left(|z|^2+|z^\prime|^2\right)\;.
 \e*
We deduce that $\Psi_{\alpha}$ attains its maximum at some
$(t_\alpha,z_\alpha,z^{\prime}_\alpha)$ in $[0,T]\x\R_+^{2d}$.
The inequality $\Psi_{\alpha}(t_0,z_0,z_0)$ $\le$
$\Psi_{\alpha}(t_\alpha,z_\alpha,z^\prime_\alpha)$ reads
 \b*
    \Psi_{\alpha}(t_\alpha,z_\alpha,z^\prime_\alpha)\;\ge\;
    \delta - \eps\left(1+4|z_0|^2\right)e^{L T}\;.
 \e*
Hence, $\eps$ can be chosen sufficiently small (depending on $L$ and $|z_0|$) so that
 \begin{eqnarray}
 \label{eq proof thm comparison 3}
 v(t_\alpha,z_\alpha)-u(t_\alpha,z^\prime_\alpha)\;\ge\;
 \Psi_{\alpha}(t_\alpha,z_\alpha,z^\prime_\alpha)\;\ge\;
 \delta - \eps\left(1+4|z_0|^2\right)e^{L T}\;>\;0\;.
 \end{eqnarray}
From \reff{eq proof thm comparison 3} and \reff{eq proof thm
comparison 2}, we get
 \b*
 0\;\le\;\frac{\alpha}{2}|z_\alpha-z^\prime_\alpha|^2+
 \frac{\eps}{2}|z_\alpha+z^\prime_\alpha|^2 &\le&
 v(t_\alpha,z_\alpha)-u(t_\alpha,z^\prime_\alpha)
 -\frac{\eps}{2} |z_\alpha+z^\prime_\alpha|^2 - \frac{\alpha}{2}|z_\alpha-z^\prime_\alpha|^2
 \\
 &\le&
 A\left( 1 + |z_\alpha| +|z^\prime_\alpha|\right) -
 \min\left\{\frac{\eps}{2},\frac{\alpha}{2}\right\}\left(|z_\alpha|^2+|z^\prime_\alpha|^2
 \right)\;.
 \e*
We deduce that
$\left\{\alpha|z_\alpha-z^\prime_\alpha|\right\}_{\alpha > 0}$ as
well as $\left\{(z_\alpha,z^\prime_\alpha)\right\}_{\alpha >0}$
are bounded. Therefore, after possibly passing to a subsequence, we
can find $(\bar t,\bar z)$ $\in$ $[0,T]\x\R_+^d$ such that
 \b*
 (t_\alpha,z_\alpha,z^\prime_\alpha) \rightarrow (\bar t,\bar
 z,\bar z) &\mbox{as}~~\alpha
 \rightarrow \infty\;.
 \e*
Since $v\,-\,u$ is upper semicontinuous, it follows from \reff{eq
proof thm comparison 3} that
 \b*
 v(\bar t,\bar z)-u(\bar t, \bar z) \;\ge\; \limsup_{\alpha \rightarrow \infty}
 v(t_\alpha,z_\alpha) - u(t_\alpha,z^\prime_\alpha) \;\ge\; \delta
 - \eps \left(1 + 4 |z_0|^2\right)e^{L T}\;>\;0\;.
 \e*
 Since $u(T,.)$ $\ge$ $v(T,.)$ on $\R_+^d$, $\bar t$ is in
 $[0,T)$, hence for $\alpha$ sufficiently large $t_\alpha$ is in
 $[0,T)$.

 \vspace{2mm}

\no 3. Let $\alpha$ be sufficiently large so that
 \b*
 |z_\alpha - z^\prime_\alpha| \;<\; 1
 &\mbox{and}& t_\alpha \,\in\, [0,T)\;.
 \e*
Since $(t_\alpha,z_ \alpha,z^\prime_\alpha)$ is a maximum point of
$\Psi_{\alpha}$, by the fundamental result in the User's Guide to
Viscosity Solutions (Theorem 8.3 in Crandall et al. 1993), for
each $\eta$ $>$ $0$, there are numbers $a_1^\eta$ , $a_2^\eta$ in
$\R$, and symmetric matrices $X^\eta$  and $Y^\eta$ in $\M^d$
such that
 \b*
 \left(a_1^\eta, D_{z}\Phi_\alpha(t_\alpha,z_\alpha,z^\prime_\alpha),
 X^\eta\right)\;\in\; \bar \Pc^{2,+}(v)(t_\alpha,z_\alpha)\;,\\
\left(a_2^\eta, -
D_{z^\prime}\Phi_\alpha(t_\alpha,z_\alpha,z^\prime_\alpha),
 Y^\eta\right)\;\in\; \bar \Pc^{2,+}(u)(t_\alpha,z^\prime_\alpha)\;,
 \e*
 with
 \b*
    a_1^\eta-a_2^\eta\;=\;\frac{\partial \Phi_{\alpha}}{\partial
    t}(t_\alpha,z_\alpha,z^\prime_\alpha)\;=\;
    - L\,\Phi_{\alpha}(t_\alpha,z_\alpha,z^\prime_\alpha)\;,
 \e*
 and
 \begin{eqnarray}
 \label{eq proof thm comparison 4}
    \left(
    \begin{array}{c c}
        X^\eta&0\\0&-Y^\eta
    \end{array}
    \right)
    \;\le\; M + \eta M^2\;,&\mbox{where}&
    M\;:=\; D^2_{(z, z^\prime)}\Phi_{\alpha}(t_\alpha,z_\alpha,z^\prime_\alpha)\;.
 \end{eqnarray}
Since $v$ is a viscosity subsolution  and $u$ is a viscosity
supersolution of \reff{eq proof thm comparison 1} on
$[0,T)\x\R_+^d$, we must have
 \b*
 \lambda v(t_\alpha,z_\alpha) - a_1^\eta - \frac{1}{2}\sup_{\mu \in
 U_\kappa}Tr\left[\sigma^\mu(t_\alpha,z_\alpha)^\prime X^\eta
 \sigma^\mu(t_\alpha,z_\alpha)\right]\;\le\;0\;,\\
\lambda u(t_\alpha,z_\alpha) - a_2^\eta - \frac{1}{2}\sup_{\mu \in
 U_\kappa}Tr\left[\sigma^\mu(t_\alpha,z^\prime_\alpha)^\prime
 Y^\eta
 \sigma^\mu(t_\alpha,z^\prime_\alpha)\right]\;\ge\;0\;.
 \e*
Taking the difference we get
 \begin{eqnarray}
 \nonumber
 &&\lambda \left( v(t_\alpha,z_\alpha)-u(t_\alpha,z^\prime_\alpha)
 \right) + L \Phi_{\alpha}(t_\alpha,z_\alpha,z^\prime_\alpha)
 \\
 \nonumber
 &\le&\frac{1}{2}\sup_{\mu \in
 U_\kappa}Tr\left[\sigma^\mu(t_\alpha,z_\alpha)^\prime X^\eta
 \sigma^\mu(t_\alpha,z_\alpha)\right]-\frac{1}{2}\sup_{\mu \in
 U_\kappa}Tr\left[\sigma^\mu(t_\alpha,z^\prime_\alpha)^\prime
 Y^\eta
 \sigma^\mu(t_\alpha,z^\prime_\alpha)\right]\\
 \label{eq proof thm comparison 5}
 &\le&
 \frac{1}{2}\sup_{\mu \in
 U_\kappa}\left\{Tr\left[\sigma^\mu(t_\alpha,z_\alpha)^\prime
 X^\eta
 \sigma^\mu(t_\alpha,z_\alpha)\right] -Tr\left[\sigma^\mu(t_\alpha,z^\prime_\alpha)^\prime
 Y^\eta
 \sigma^\mu(t_\alpha,z^\prime_\alpha)\right]\right\}\;.
 \end{eqnarray}
Let $(e_i, i=1,...,d)$ be an orthonormal basis of $\R^d$, and for
each $\mu$ in $U_\kappa$ set
 \b*
 \xi^\mu_i\;:=\;\sigma^\mu(t_\alpha,z_\alpha)e_i&\mbox{and}&\zeta^\mu_i\;:=\;\sigma
 ^\mu(t_\alpha,z^\prime_\alpha)e_i
 \e*
 so that
 \b*
&&Tr\left[\sigma^\mu(t_\alpha,z_\alpha)^\prime
X^\eta\sigma^\mu(t_\alpha,z_\alpha)\right]
    -Tr\left[\sigma^\mu(t_\alpha,z^\prime_\alpha)^\prime Y^\eta
    \sigma^\mu(t_\alpha,z^\prime_\alpha)\right]
 \\&=&\sum_{i=1}^d
 X^\eta\xi^\mu_i \cdot \xi^\mu_i-
 Y^\eta\zeta^\mu_i\cdot \zeta^\mu_i
 \e*
and, by \reff{eq proof thm comparison 4},
 \b*
Tr\left[\sigma^\mu(t_\alpha,z_\alpha)^\prime
X^\eta\sigma^\mu(t_\alpha,z_\alpha)\right]
    -Tr\left[\sigma^\mu(t_\alpha,z^\prime_\alpha)^\prime Y^\eta
   \sigma^\mu(t_\alpha,z^\prime_\alpha)\right]
 &\le&
 \sum_{i=1}^d (M+\eta M^2) \beta^\mu_i
 \cdot \beta^\mu_i \;,
 \e*
where $\beta^\mu_i$ is the $2d$-dimensional column vector
defined by ~: ~$\beta^\mu_i$ $:=$
$(\xi^{\mu\,\prime}_i,\zeta^{\mu \, \prime}_i)^\prime$.
  Letting $\eta$ go to zero, and using \reff{eq proof thm
comparison 5}, we get
 \begin{equation}
 \label{eq proof thm comparison 6}
 \lambda \left( v(t_\alpha,z_\alpha)-u(t_\alpha,z^\prime_\alpha)
 \right) + L \Phi_{\alpha}(t_\alpha,z_\alpha,z^\prime_\alpha)
 \;\le\;\frac{1}{2}\sup_{\mu \in
 U_\kappa}\left\{
 \sum_{i=1}^d M\beta^\mu_i
 \cdot\beta^\mu_i
 \right\}\;.
 \end{equation}

 \vspace{2mm}

 \no 4. In this last step, we are going to see  that, for a
 convenient choice of the positive constant $L$, inequality \reff{eq proof thm comparison
 6} leads to a contradiction to \reff{eq proof thm comparison 3}.

 \no Notice that
 \b*
 M\;=\;e^{L (T-t_\alpha)}\left(
 \begin{array}{l l}
 D_{zz}\phi_{\alpha}+D_{z^\prime
 z^\prime}\phi_{\alpha}+2 D_{z z^\prime} \phi_{\alpha}&
 D_{z z} \phi_{\alpha}-D_{z^\prime z^\prime}\phi_{\alpha}\\
 D_{z z} \phi_{\alpha}-D_{z^\prime z^\prime}\phi_{\alpha}&
 D_{zz}\phi_{\alpha}+D_{z^\prime
 z^\prime}\phi_{\alpha}-2 D_{z z^\prime} \phi_{\alpha}
 \end{array}
 \right)
 \e*
the (partial) Hessian matrices of $\phi_{\alpha}$ being taken at
the point $(z_\alpha+z^\prime_\alpha,z_\alpha -z^\prime_\alpha)$.
Then, for $\mu$ in $U_\kappa$
 \begin{eqnarray}
\nonumber
 \sum_{i=1}^d M\beta^\mu_i
 \cdot \beta^\mu_i
 &=&
 e^{L(T-t_\alpha)}\sum_{i=1}^d D_{z
 z}\phi_{\alpha}\left(\xi^\mu_i+\zeta^\mu_i\right)\cdot\left(\xi^\mu_i+\zeta^\mu_i\right)
 \\
 \label{eq proof thm comparison 10}
 &&+\; 2D_{z
 z^\prime}\phi_{\alpha}\left(\xi^\mu_i+\zeta^\mu_i\right)\cdot\left(\xi^\mu_i-\zeta^\mu_i\right)\\
\nonumber
 &&+\;  D_{z^\prime
 z^\prime}\phi_{\alpha}\left(\xi^\mu_i-\zeta^\mu_i\right)\cdot\left(\xi^\mu_i-\zeta^\mu_i\right)\;.
 \end{eqnarray}

\vspace{2mm}
 \no Since for each $\mu$ in $U_\kappa$, $|\mu|$ is
bounded by $\kappa$, we verify using the Lipschitz property of the
function $z$ $\mapsto$ $diag[z]\sigma(t,z)$ that for some
positive constant $C$, for each $z$, $z^\prime$ in $\R_+^d$ and
$t$ in $[0,T]$
 \b*
|\sigma^\mu(t,z) - \sigma^\mu(t,z^\prime)|\;\le\; C|z- z^\prime|
&\mbox{and}& |\sigma^\mu(t,z)|\;\le\; C\left(1+|z|\right)\;.
 \e*

\vspace{2mm}

\no In the following $C$ denotes a positive constant
(independent of $\alpha$, $\eps$ and $L$) which value may change
from line to line.

\no 4.1. Since $\alpha$ satisfies $|z_\alpha -
 z^\prime_\alpha|$ $\le$ $1$,  for $i$ $=$ $1,..,d$,
 \b*
 D_{z
 z}\phi_{\alpha}\left(\xi^\mu_i+\zeta^\mu_i\right)\cdot\left(\xi^\mu_i+\zeta^\mu_i\right)
 \;\le\;
 C |D_{z
 z}\phi_{\alpha}|\left[(1+|z_\alpha|)^2+
 (1+|z^\prime_\alpha|)^2\right]\;.
 \e*
From $|D_{z
 z}\phi_{\alpha}|$ $\le$ $2(\eps + \alpha|z_\alpha
 -z^\prime_\alpha|^2)$ and the previous estimate, we deduce
 that
\begin{eqnarray}
\label{eq proof thm comparison 7}
 D_{z
 z}\phi_{\alpha}\left(\xi^\mu_i+\zeta^\mu_i\right)\cdot\left(\xi^\mu_i+\zeta^\mu_i\right)
 &\le&
 C(\eps + \alpha|z_\alpha
 -z^\prime_\alpha|^2)\left(1 + |z_\alpha +
 z^\prime_\alpha|^2\right)\;.
 \end{eqnarray}

\no 4.2. For $i$ $=$ $1,..,d$
 \b*
D_{z
 z^\prime}\phi_{\alpha}\left(\xi^\mu_i+\zeta^\mu_i\right)\cdot\left(\xi^\mu_i-\zeta^\mu_i\right)
 \;\le\;
 C |D_{z
 z^\prime}\phi_{\alpha}|\left[(1 + |z_\alpha|) + (1+
 |z^\prime_\alpha| )\right]|z_\alpha - z^\prime_\alpha|\;.
 \e*

\no Since $|D_{z
 z^\prime}\phi_{\alpha}|$ $\le$ $4\alpha|z_\alpha -z^\prime_\alpha||z_\alpha
 +z^\prime_\alpha|$ and $|z_\alpha -
 z^\prime_\alpha|$ $\le$ $1$, we deduce that
 \begin{eqnarray}
 \label{eq proof thm comparison 8}
D_{z
 z^\prime}\phi_{\alpha}\left(\xi^\mu_i+\zeta^\mu_i\right)\cdot\left(\xi^\mu_i-\zeta^\mu_i\right)
 &\le&
 C(\eps + \alpha|z_\alpha
 -z^\prime_\alpha|^2)\left(1 + |z_\alpha +
 z^\prime_\alpha|^2 \right)\;.
 \end{eqnarray}

 \no 4.3. For $i$ $=$ $1,..,d$
 \b*
D_{z^\prime
z^\prime}\phi_{\alpha}\left(\xi^\mu_i-\zeta^\mu_i\right)\cdot\left(\xi^\mu_i-\zeta^\mu_i\right)
\;\le\; C|D_{z^\prime z^\prime}\phi_{\alpha}||z_\alpha-
z^\prime_\alpha|^2
 \e*

\no and since $|D_{z^\prime z^\prime}\phi_{\alpha}|$ $\le$
$2\alpha\left(1+ |z_\alpha+z^\prime_\alpha|^2\right)$, we get
 \begin{eqnarray}
 \label{eq proof thm comparison 9}
 D_{z^\prime
z^\prime}\phi_{\alpha}\left(\xi^\mu_i-\zeta^\mu_i\right)\cdot\left(\xi^\mu_i-\zeta^\mu_i\right)
&\le& C (\eps + \alpha|z_\alpha
 -z^\prime_\alpha|^2)\left(1 + |z_\alpha +
 z^\prime_\alpha|^2 \right)\;.
 \end{eqnarray}

\no Finally, collecting the estimates \reff{eq proof thm
comparison 7}, \reff{eq proof thm comparison 8}  and \reff{eq
proof thm comparison 9}, we deduce from \reff{eq proof thm
comparison 10} that for some positive constant $\tilde C$
(independent of $L$, $\eps$ and $\alpha$)
 \b*
 \sum_{i=1}^d M \beta^\mu_i
 \cdot \beta^\mu_i
 &\le& \tilde C e^{L(t-t_\alpha)}(\eps + \alpha|z_\alpha
 -z^\prime_\alpha|^2)\left(1 + |z_\alpha +
 z^\prime_\alpha|^2 \right)\;=\; \tilde C
 \Phi_{\alpha}(t_\alpha,z_\alpha,z^\prime_\alpha)\;.
 \e*

\no Hence, if we take $L$ $\ge$ $\frac{\tilde C}{2}$, then
\reff{eq proof thm comparison 6} reads
 \b*
 \lambda (v(t_\alpha,z_\alpha) - u(t_\alpha,z^\prime_\alpha)) \;\le\; (\frac{\tilde C}{2} - L) \Phi_{\alpha}(t_\alpha,z_\alpha,z^\prime_\alpha) \;\le\;0
 \e*
 which is in contradiction with  \reff{eq proof thm comparison 3}.

 \ep

\subsection{ Proof of Proposition \ref{prop u ge w}}

\no We first make use of Proposition \ref{prop theo comp} to
obtain the intermediary inequality $u_*\ge \sup_{\kappa>0}
v^*_\kappa$.

\begin{Corollary}\label{cor u ge vkappa} Under $(H\lambda)$, $(H\sigma)$ and $(Hg)$, for each $\kappa>0$, we
have $u_*\ge v^*_\kappa$ on $[0,T] \times \R_+^d$.
\end{Corollary}

\proof In view of Lemmas \ref{lem vkappa soussol}, \ref{lem
vkappa soussol cond bord}, Corollary \ref{cor BT} and Proposition
\ref{prop theo comp}, it suffices to show that $u_*$ and
$v^*_\kappa$ have linear growth.  To check this condition for
$v^*_\kappa$, it suffices to recall that  $Z^\mu$ is a martingale
and use assumption $(Hg)$. We now consider $u_*$. First, recall
from $(Hg)$ that $\hat G$ has linear growth. Using Lemma \ref{lem
f}, we deduce that, for each $(a,b) \in \Dc$ and $(t,y,z) \in
[0,T]\x (0,\infty)^{n+d}$, we have
    \b*
    f\left(Y^{(a,b)}_{t,y,z}(T)\right)\cdot g\left(S_{t,s}(T)\right)
    &\le&
    \hat G\left(Z^{(a,b)}_{t,y,z}(T)\right)
    \;\le\;
    \delta \left(1+  \sum_{i=1}^d Z^{(a,b)^i}_{t,y,z}(T) \right)
    \e*
where $s$ $=$ $\bar F(y)^{-1}z$ and $\delta$ is some positive constant. It follows from the definition of $u^{(a,b)}$ that
    \be\label{eq borne sup u}
    u(t,y,z)\;=\;   \sup_{(a,b) \in \Dc} u^{(a,b)}(t,y,z) \le  \delta   \left(1+  \sum_{i=1}^d z^i\right)\;.
    \ee
Now, observe that, for $(a,b)=(0,0)$, $Y^{(0,0)}$ is constant so
that $Z^{(0,0)^i}$ coincides with $S^i$ up to a multiplicative
constant, $i\le d$. Hence, $Z^{(0,0)}$ is a $\P$-martingale and
it follows from the definition of $u^{(0,0)}$ that
    \b*
    u(t,y,z)
    &\ge&
    u^{(0,0)}(t,y,z)
    \;\ge\;
    \Esp{   f\left(y\right)\cdot
    g\left(S_{t,s}(T)\right)}\;\;,\;s=\bar F(y)^{-1}z\;.
    \e*
Using  Lemma \ref{lem f} and \reff{TA}, we then deduce as above  that
    \b*
    u(t,y,z)
    &\ge&
    u^{(0,0)}(t,y,z)
    \;\ge\;
    -\hat \delta \left(1+ \sum_{i=1}^dz^i\right)
    \e*
for some positive constant  $\hat \delta$. Combining the last inequality with \reff{eq borne sup u} shows that $u_*$ has linear growth.
 \ep
\\

\no{\bf Proof of Proposition \ref{prop u ge w}.} 1. We first show
that $\{Z_{0,z}^{\mu}\;,\;\mu \in \cup_{\kappa\ge 0}
\Uc_\kappa\}$ is dense in probability in $\{Z_{0,z}^{\mu}\;,\;\mu
\in   \Uc \}$. To see this,  take $\mu$ $\in \Uc$ and consider
the sequence defined by $\mu_\kappa:=\mu \1_{|\mu| \le \kappa}
\in \Uc_\kappa$, $\kappa \in \N$. Recalling  that
$Z^{\mu,i}_{0,z}$ $=$ $S^i_{0,z}$ for $i\le d_f$, see assumption
$(H\sigma)$, we deduce from the It\^o's isometry that
    \b*
    &&\Esp{\sum_{i=1}^d |\ln (Z_{0,z}^{\mu,i}(T))-\ln (Z_{0,z}^{\mu_\kappa,i}(T))|^2}
    \\
    &&\le\;
    \delta \Esp{\int_0^T |\mu(t)-\mu_\kappa(t)|^2 + ||\mu(t)|^2-|\mu_\kappa(t)|^2| \;dt}\;
    \e*
for some $\delta>0$.
Since $\mu_\kappa \to \mu$ ${\rm dt}\times {\rm d}\P$-a.e. and, by definition of $\Uc$, $\mu$ is square integrable, we deduce from
the dominated convergence theorem that  $\ln (Z_{0,z}^{\mu_\kappa,i}(T))$
goes to $\ln  (Z_{0,z}^{\mu,i}(T))$ in $L^2$, $i\le d$. It follows that the convergence holds $\Pas$   along some
subsequence.

\no 2.  By  Corollary \ref{cor u ge vkappa}, we have
    \b*
    u_*(0,z)&\ge&  \sup_{\kappa>0} \sup_{\mu \in \Uc_\kappa} \Esp{\hat  G(Z_{0,z}^\mu(T))}\;.
    \e*
Since $\hat G$ has linear growth, see $(Hg)$, there is some
$(c,\Delta) \in \R\x \R^{d}$ such that $\hat  G(Z_{0,z}^\mu(T)) +
\Delta\cdot Z_{0,z}^\mu(T) \ge -c$. Since $Z_{0,z}^\mu$ is a
martingale, it follows that
      \b*
    u_*(0,z)&\ge&
    \sup_{\kappa>0}
    \sup_{\mu \in \Uc_\kappa}
    \Esp{\hat  G(Z_{0,z}^\mu(T))+ \Delta\cdot Z_{0,z}^\mu(T)}
    - \Delta\cdot z\;.
    \e*
Using 1. and Fatou's Lemma, we then deduce that
      \b*
    u_*(0,z)&\ge&
    \sup_{\mu \in \Uc}
    \Esp{\hat  G(Z_{0,z}^\mu(T))+ \Delta\cdot Z_{0,z}^\mu(T)}
    - \Delta\cdot z\;.
    \e*
Since for $\mu \in \Uc$, $Z_{0,z}^\mu$ is also a martingale, the result follows.   \ep

\section{Proof of Theorem \ref{thmmain} and Corollary \ref{cor formulation duale}}\label{sec proof finale thmmain}

\no In order to conclude the proof,  we shall now exploit the
concavity property of Corollary \ref{cor BT}.

\begin{Remark}{\rm Fix $\mu \in \Uc$, $t\le T$ and $(s^f,z^c)\in (0,\infty)^{d_f} \x
(0,\infty)^{d_c}$. Then, it follows from assumption $(H\sigma)$
that $Z^{\mu,f}_{t,(s^f,z^c)}=S^f_{t,(s^f,z^c)}$ which depends only on
$(t,s^f)$. On the other hand $Z^{\mu,c}_{t,(s^f,z^c)}$ depends
only on $(t,z^c,\mu)$. In the following, we shall then slightly
abuse notations and simply write $S^f_{t,s^f}$ for $S^f_{t,s}$ and
$Z^{\mu,c}_{t,z^c}$ for $Z^{\mu,c}_{t,(s^f,z^c)}$. }
\end{Remark}

\begin{Corollary}\label{cor borne finale pour u} Let the conditions $(H\lambda)$, $(H\sigma)$ and $(Hg)$ hold.
Then, for all $(s^f,z^c)\in (0,\infty)^{d_f} \x (0,\infty)^{d_c}$
and $\Delta\in \partial_{z^c}u_*(0,s^f,z^c)$, we have
    \b*
    u_*(0,s^f,z^c)
    &\ge&
    \Esp{
    \sup_{\tilde z^c \in (0,\infty)^{d_c}}
    \left\{\hat  G(S^f_{0,s^f}(T),\tilde z^c)-\Delta\cdot  \tilde z^c\right\}
    }
    +\Delta\cdot  z^c\;.
    \e*
\end{Corollary}

\proof By definition of $\Delta$ and Corollary \ref{cor BT}, we
have
    \b*
    u_*(0,s^f,z^c) &=& \sup_{\tilde z^c\in (0,\infty)^{d_c}}\;\left\{ u_*(0,s^f,\tilde z^c)
                                            - \Delta\cdot  (\tilde  z^c - z^c)\right\}\;.
    \e*
Since, for each $\tilde z^c\in (0,\infty)^{d_c}$ and $\mu \in \Uc$,
$\Esp{Z^{\mu,c}_{0,\tilde z^c}}=\tilde z^c$, it follows from
Proposition \ref{prop u ge w} that
    \b*
    u_*(0,s^f,z^c)
    &\ge&
    \sup_{\tilde z^c\in (0,\infty)^{d_c}}\sup_{\mu \in \Uc}
    \Esp{\hat  G\left(S^f_{0,s^f}(T),Z_{0,\tilde z^c}^{\mu,c}(T)\right)-\Delta\cdot Z_{0,\tilde z^c}^{\mu,c}(T)}
    +\Delta\cdot  z^c\;.
    \e*
Since $\hat G$ is continuous, we deduce from the representation
theorem  that
    \b*
    u_*(0,s^f,z^c)
    &\ge&
    \sup_{ \xi \in L^\infty(\Bc_\kappa;\Fc_T)}
    \Esp{\hat  G(S^f_{0,s^f}(T),\xi)-\Delta\cdot \xi}
    +\Delta\cdot  z^c\;
    \\
    &\ge &
    \Esp{\max_{\tilde z^c\in \Bc_\kappa}
    \hat  G(S^f_{0,s^f}(T),\tilde z^c)-\Delta\cdot \tilde z^c}
    +\Delta\cdot  z^c\;,
    \e*
where $\Bc_\kappa$ $:=$ $\{\alpha \in (0,\infty)^{d_c}~:~$
$|\ln(\alpha^i)| \le \kappa\;\;,\;i\le d_c\}$.
 The result then follows from monotone convergence.
 \ep
\\

\no We can now conclude the proof of our main result. \\

\no{\bf Proof of Theorem \ref{thmmain}.}  For ease of notations, we
only write $S$ for $S_{0,S(0)}$ and $s=(s^f,s^c)$ for $S(0)$. In view of Theorem \ref{thm BT}, we have
    \be\label{eq *}
    p(0,S(0))
    &\ge&
    \sup_{\xi\in \Lambda}
    u_*(0,s^f,\diag[\lb{\xi}]s^c)\;.
    \ee

\no 1. Recalling that $\Lambda$ is compact and $u_*$ is concave in
its last $d_c$ variables, there is some $\hat \xi$ $\in \Lambda$
which attains the optimum in the above inequality. Moreover, by standards arguments of calculus of variations,
we can find some $\hat \Delta \in \partial_{z^c}u_*(0,s^f,\diag[\lb{\hat \xi}]s^c)$ such that
    \be\label{eq achat init finance}
    (\diag[\hat \Delta]s^c)\cdot ( \lb{\hat \xi}  -   \lb{\xi} )&\ge& 0   \mbox{
for all } \xi \in \Lambda\;.
    \ee
From   \reff{eq cara K par Lambda}, we deduce that $(\diag[\hat
\Delta]s^c\cdot  \lb{\hat \xi},0)\succeq (0,\diag[\hat
\Delta]s^c)$.

\no 2. Set
        \b*
        \delta\;:=\; \diag[\hat \Delta]s^c\cdot  \lb{\hat \xi} &\And&
        \hat C(S^f(T))\;:=\;\sup_{\tilde z \in (0,\infty)^{d_c}} \hat  G(S^f(T),\tilde z)-\hat \Delta\cdot  \tilde
            z \;\;,
       \e*
so that, by \reff{eq *} and Corollary \ref{cor borne finale pour u},
   \be\label{eq **}
    p(0,S(0))&\ge& p\;:=\;
    \E\left[\hat C(S^f(T))\right]+\delta\;.
    \ee
Since by $(H\sigma)$ the dynamics of $S^f$ depends only on
$S^f$,  it follows that there is some $\R^{d_f}$-valued
predictable process $\phi$ satisfying $\int_0^T |\phi(t)|^2 dt
<\infty$  such that
    \begin{equation}\label{eq def Xf a la fin}
    X^f(\cdot)\;:=\;p-\delta
    +
    \int_0^\cdot \phi(t)\cdot  dS^f(t)
    \mbox{ is a  martingale and }
     X^f(T)
    \;=\;
    \hat C(S^f(T)) \;.
    \end{equation}

\no 3. By combining 1. and 2., we  deduce that there is some strategy $(\phi,L)$ such that $L(t)$ $=$ $L(0)$,  $X^{\phi,L}_{p\1_1}(0)$ $=$ $(p-\delta,\diag[\hat \Delta]s^c)$ and
    \b*
    X^{\phi,L}_{p\1_1}(T)
    &=&
    \left(X^f(T),\diag[\hat \Delta]S^c(T)\right)
    \;\succeq\;
    \left(\hat C(S^f(T)),\diag[\hat \Delta]S^c(T)\right) \;.\;
    \e*
Using  \reff{eq def Xf a la fin}, \reff{eq cara K par Lambda} and
the definition of $\hat C$ this implies that
    \b*
    \tilde \xi\cdot X^{\phi,L}_{p\1_1}(T) -\xi\cdot  g\left(S^f(T),\diag[ \lb{\xi}]^{-1}\diag[\lb{\tilde \xi}]S^c(T)\right)
    &\ge& 0\;\;\pourtout\;\tilde \xi,\;\xi \in \Lambda\;\;.\;
    \e*
Considering the case where $\xi$ $=$ $\tilde \xi$ and using
\reff{eq cara K par Lambda} leads to
        \b*
    X^{\phi,L}_{p\1_1}(T) &\succeq& g(S(T)) \;.
    \e*
In view of \reff{eq **} it remains to check that $(\phi,L) \in \Ac$, but this readily follows from \reff{eq def Xf a la fin}
and assumption \reff{TA}.
 \ep
\\

\no {\bf Proof of Corollary \ref{cor formulation duale}.}  In
view of  Theorem \ref{thmmain}, we only have to show that
    \b*
    p(0,S(0))&\le&
    \inf_{\Delta \in \R^{d_c}}
    \;
    \Esp{C(S^f(T);\Delta)}
    + \sup_{\xi \in \Lambda} \lb{\xi}\cdot  \diag[ \Delta ] S^c(0)
    \;.
    \e*
To see this, fix some $\Delta \in \R^{d_c}$ such that
    \b*
    \tilde p &:= & \Esp{C(S^f(T);\Delta)}
    + \sup_{\xi \in \Lambda} \lb{\xi}\cdot  \diag[ \Delta ] S^c(0)
    \;<\; \infty\;,
    \e*
which is possible by  Theorem \ref{thmmain}. Then, by the same
argument as in the proof of Theorem \ref{thmmain}, see above, we
obtain that there exists some $(\phi, L)$ $\in \Ac^{BH}$, with
$L(t)=\Delta$ on $t\le T$, such that $X^{\phi,L}_{\tilde p\1_1}
\succeq g(S(T))$. This proves the required inequality as well as
the last statement of the Corollary. \ep


\section*{References}

\no {\sc Barles, B., Biton, S., Bourgoing, M. and Ley ,O. } (2003),
``Quasilinear parabolic equations, unbounded solutions and
geometrical equations III. Uniqueness through classical viscosity
solution's methods ", Calc. Var. Partial Differential Equations,
{ 18}, 159-179.

\vspace{3mm}

\no {\sc Bouchard, B. } (2000), {\sl Contr\^ole stochastique appliqu\'e \`a la finance}, Phd Thesis, University Paris-Dauphine.

\vspace{3mm}

\no {\sc Bouchard, B. and Touzi, N.} (2000), ``Explicit solution
of the multivariate super-replication problem under transaction
costs", The Annals  of Applied Probability, { 10} 3, 685-708.

\vspace{3mm}

\noindent {\sc Crandall, M.G., Ishii, H. and Lions, P.L.} (1992), ``User's
guide to viscosity solutions of second order Partial Differential Equations",
Bull. Amer. Math. Soc. { 27}, 1-67.

\vspace{3mm}

\noindent {\sc Cvitani\'c, J. and Karatzas, I.} (1996), ``Hedging and
portfolio optimization under transaction costs", {\sl Mathematical Finance},
{ 6}, 133-165.

\vspace{3mm}

\noindent {\sc Cvitani\'c, J., Pham, H. and Touzi, N.} (1999), ``A closed-form
solution to the problem of super-replication under transaction costs",
{\sl Finance and Stochastics}, { 3}, 35-54.

\vspace{3mm}

\noindent {\sc Davis, M. and Clark, J.M.C.} (1994), ``A note on
super-replicating
strategies", {\sl Phil. Trans. Roy. Soc. London A}, { 347}, 485-494.

\vspace{3mm}

\noindent {\sc Jouini, E. and Kallal, H.} (1995), ``Martingales and Arbitrage
in securities markets with transaction costs", {\sl Journal of Economic
Theory},
{ 66}, 178-197.

\vspace{3mm}

\noindent {\sc Kabanov, Yu.} (1999), ``Hedging and liquidation under
transaction
costs in currency markets", {\sl Finance and Stochastics}, { 3} (2), 237-248.

\vspace{3mm}

\noindent {\sc Levental, S. and Skorohod, A.V.} (1997), ``On the possibility of
hedging options in the presence of transaction costs", {\sl Annals of Applied
Probability}, { 7}, 410-443.

\vspace{3mm}

\noindent {\sc Rockafellar, R.T.} (1970), {\sl Convex Analysis}, Princeton
University Press, Princeton, NJ.

\vspace{3mm}

\noindent {\sc Soner, M., Shreve, S. and Cvitani\'c J.} (1995), ``There is no
nontrivial hedging portfolio for option pricing with transaction costs",
{\sl Annals of Applied Probability}, { 5}, 327-355.

\vspace{3mm}

\noindent {\sc Soner, M.  and Touzi, N.} (2002), ``Dynamic programming for stochastic target problems and geometric flows",
{\sl Journal of the European Mathematical Society},   4, 201-236.

\vspace{3mm}

\noindent {\sc Touzi, N.}  (1999), ``Super-replication under
proportional transaction costs: from discrete to continuous-time
models", {\sl Mathematical Methods of Operations Research}, 50,
297-320.

\end{document}